\documentstyle{article}
\font\germ=eufm10
\def\ssl{\hbox{\germ sl}}
\def\slh{\widehat{\ssl_2}}
\makeatletter
\def\aaa{@}
\title{\Large\bf Crystallized Peter-Weyl type decomposition
for level 0 part of modified quantum algebra
$\widetilde U_q(\slh)_0$}
\author{Toshiki Nakashima
\thanks{supported by
the Ministry of Education, Science and Culture of Japan
as a overseas research scholar.
Current address : Department of Mathematics,
Northeastern University, Boston, MA 02115, USA.
(e-mail : toshiki\aaa neu.edu)}
\\
Department of Mathematical Science, Faculty of Engineering Science,\\
Osaka University, Toyonaka, Osaka 560, Japan\\
toshiki\aaa sigmath.es.osaka-u.ac.jp}
\date{}
\makeatother
\begin{document}

\maketitle

   \renewcommand{\labelenumi}{$(\roman{enumi})$}
  \font\germ=eufm10

  \def\aff{{\rm Aff}}
  \def\al{\alpha}
  \def\beq{\begin{equation}}
  \def\beqn{\begin{eqnarray}}
  \def\beqnn{\begin{eqnarray*}}
  \def\bigsl{{\hbox{\fontD \char'54}}}
  \def\binf{B_{\infty}}
  \def\bminf{B_{-\infty}}
  \def\bmax{B^{\scriptstyle{\rm max}}}
  \def\bmin{B^{\scriptstyle{\rm min}}}
  \def\catob{{\cal O}(B)}
  \def\cd{\cdots}
  \def\del{\delta}
  \def\Del{\Delta}
  \def\ei{e_i}
  \def\eit{\tilde{e}_i}
  \def\ep{\epsilon}
  \def\eeq{\end{equation}}
  \def\eeqn{\end{eqnarray}}
  \def\eeqnn{\end{eqnarray*}}
  \def\FF{\hbox{\bf F}}
  \def\fit{\tilde{f}_i}
  \def\ge{\hbox{\germ g}}
  \def\gl{\hbox{\germ gl}}
  \def\gc{\bigcirc}
  \def\hom{{\hbox{Hom}}}
  \def\ify{\infty}
  \def\io{\iota}
  \def\llra{\relbar\joinrel\relbar\joinrel\relbar\joinrel\rightarrow}
  \def\lan{\langle}
  \def\lar{\longrightarrow}
  \def\lm{\lambda}
  \def\Lm{\Lambda}
  \def\mapright#1{\smash{\mathop{\longrightarrow}\limits^{#1}}}
  \def\mapleftright#1{\smash{\mathop{\longleftrightarrow}\limits^{#1}}}
  \def\map#1{\smash{\mathop{\longmapsto}\limits^{#1}}}
  \def\mpath{{\cal P}_m}
  \def\mlpath{{\cal P}_{m,l}}
  \def\mpathtc{{\cal P}_m(n;\vec{t};\vec{c})}
  \def\mlpathtc{{\cal P}_{m,l}(n;\vec{t};\vec{c})}
  \def\mllpathtc{{\cal P}_{m,l'}(n;\vec{t};\vec{c})}
  \def\nn{\nonumber}
  \def\nd{\noindent}
  \def\ot{\otimes}
  \def\op{\oplus}
  \def\ovl{\overline}
  \def\qq{\qquad}
  \def\q{\quad}
  \def\qed{\hfill\framebox[2mm]{}}
  \def\QQ{\hbox{\bf Q}}
  \def\qi{q_i}
  \def\qii{q_i^{-1}}
  \def\ran{\rangle}
  \def\SFF{\scriptstyle\hbox{\bf F}}
  \def\ssl{\hbox{\germ sl}}
  \def\slh{\widehat{\ssl_2}}
  \def\syl{\scriptstyle}
  \def\ti{t_i}
  \def\tii{t_i^{-1}}
  \def\til{\tilde}
  \def\tt{{\hbox{\germ{t}}}}
  \def\ttt{\hbox{\germ t}}
  \def\uq{U_q(\ge)}
  \def\uqm{U^-_q(\ge)}
  \def\uqmq{{U^-_q(\ge)}_{\bf Q}}
  \def\uqq{U^{\bf Q}_q(\ge)}
  \def\util{\widetilde\uq}
  \def\vmax{V^{\scriptstyle{\rm max}}}
  \def\vep{\varepsilon}
  \def\vp{\varphi}
  \def\vpi{\varphi^{-1}}
  \def\wtil{\widetilde}
  \def\what{\widehat}
  \def\ZZ{\hbox{\bf Z}}
%\vspace{-26pt}
%\begin{abstract}

%\end{abstract}

\section{Introduction}
\setcounter{equation}{0}
\renewcommand{\theequation}{\thesection.\arabic{equation}}

Modified quantum algebra $\util$ is an algebra which
is 'modified' the Cartan part of the underlying quantum algebra
$\uq$ as in (\ref{ualm}) and (\ref{eqn:def-util}).
Modified quantum algebra $\util$ holds the
remarkable property that
modified quantum algebra $\util$ affords the commutative
two crystal structures. The first one is the usual crystal structure,
which was found by Lusztig (\cite{L1}).
Another crystal structure was discovered by Kashiwara (\cite{K4}),
which is called right crystal structure.
In \cite{K4}, it is shown that
$\util$ is stable by the action of the antiautomorphism $*$
(see 2.1) and moreover, crystal base $B(\util)$ is also
stable by the action of $*$. By using $*$,
right crystal structure is constructed.
The commutativity of those crystal structures motivated us to
consider Peter-Weyl type decomposition on the crystal base
of $\util$.

In \cite{K4}, Kashiwara gave the Peter-Weyl type decomposition
for the crystal base of modified quantum
algebra of finite type and affine type of non-zero level part.
But the Peter-Weyl type decomposition for
affine type with level 0 part is still unclear.
In this paper, we shall give some criteria for the existence
of the Peter-Weyl type decomposition:
$$
B(\util)\cong
\bigoplus_{\lm\in P/W}
\bmax(\lm)\ot B(-\lm)^*,
$$
where $P$ is a weight lattice, $W$ is the Weyl group associated with $\ge$
and $\bmax(\lm)$ and $B(-\lm)^*$ will be given in section 4.
Those criteria are related to the property of
connected component in $B(\util)$.
Furthermore, we can consider
its application to the level 0 part of modified quantum affine algebra
$\wtil U_q(\slh)_0$ since we have already classified all
connected components in $\wtil U_q(\slh)_0$ in \cite{N} and
got the properties required by the criteria.

Let us see the organization of this paper.
In section 2, we prepare the notion of crystals and related subjects.
In section 3, we review the definition of
modified quantum algebras, the properties of their crystal base and
general feature of operation $*$ and Weyl group on $B(\util)$.
In this section, we shall give the definition of extremal vector.
In section 4,
we shall give the right structure on $B(\util)$.
Then we shall investigate the
criteria for the existence of the Peter-Weyl type decomposition.
In section 5, we consider the application  of the criteria to the
level 0 part of $B(\wtil U_q(\slh))$. In order to
give the explicit form of $\bmax(\lm)$ and $B(-\lm)^*$,
we describe the action of $*$ on extremal vectors.

The author would like to acknowledge professor Masaki Kashiwara
for his helpful advises.
This work was partly done during the stay of the author
at Northeastern University.
He is grateful to professor Andrei Zelevinsky
for his kind hospitality.

%%%%%% Section 2 %%%%%%
\section{Crystals}
\setcounter{equation}{0}
\renewcommand{\theequation}{\thesection.\arabic{equation}}
\subsection{Definition of $U_q(\ge)$}

Let $\ge$ be
a  symmetrizable Kac-Moody algebra over {\bf Q}
with a Cartan subalgebra
$\ttt$, $\{\al_i\in\ttt^*\}_{i\in I}$
 the set of simple roots and
$\{h_i\in\ttt\}_{i\in I}$  the set of coroots,
where $I$ is a finite index set. We define an inner product on
$\ttt^*$ such that $(\al_i,\al_i)\in{\bf Z}_{\geq 0}$ and
$\lan h_i,\lm\ran=2(\al_i,\lm)/(\al_i,\al_i)$
for $\lm\in\ttt^*$.
Set $Q=\oplus_i\ZZ\al_i$,
$Q_+=\oplus_i\ZZ_{\geq0}\al_i$ and
$Q_-=-Q_+$. We call $Q$ a root lattice.
Let  $P$  a lattice of $\ttt^*$ {\it i.e.} a free
{\bf Z}-submodule of $\ttt^*$ such that
$\ttt^*\cong {\hbox{\bf Q}}\ot_{\ZZ}P$,
and $P^*=\{h\in \ttt|\lan h,P\ran\subset\ZZ\}$.
We set $P_+=\{\lm\in P|\lan \lm,h_i\ran\geq 0
{\hbox{ for any }}i\in I\}$.
An element of $P$(resp.$P_+$) is called a integral weight
(resp. dominant integral weight).

The quantized enveloping algebra $\uq$ is an associative
$\QQ(q)$-algebra generated by $e_i$, $f_i(i\in I)$
and $q^h(h\in P^*)$
satisfying the following relations:
\begin{eqnarray}
&&q^0=1, \q{\hbox{\rm and }}\q q^hq^{h'}=q^{h+h'},\\
&&q^he_iq^{-h}=q^{\lan h,\al_i\ran}e_i,\qq
q^hf_iq^{-h}=q^{-\lan h,\al_i\ran}f_i,\\
&&[e_i,f_j]=\del_{i,j}(t_i-t^{-1}_i)/(q_i-q^{-1}_i),\\
&&\sum_{k=1}^{1-{\lan h_i,\al_j\ran}}
(-1)^kx_i^{(k)}x_jx_i^{(1-{\lan h_i,\al_j\ran}-k)}=0,\qq(i\ne j)
\end{eqnarray}
where $x_i=e_i,f_i$ and we set $q_i=q^{(\al_i,\al_i)/2}$,
$t_i=q_i^{h_i}$, $[n]_i=(q^n_i-q^{-n}_i)/(q_i-q_i^{-1})$,
$[n]_i!=\prod_{k=1}^n[k]_i$,
$e_i^{(n)}=e_i^n/[n]_i!$ and $f_i^{(n)}=f_i^n/[n]_i!$.

It is well-known that $\uq$ has a Hopf algebra structure with a
comultiplication $\Del$ given by
$$
\Del(q^h)=q^h\ot q^h,\q \Del(\ei)=\ei\ot\tii+1\ot\ei,\q
\Del(f_i)=f_i\ot 1+\ti\ot f_i,
$$
for any $i\in I$ and $h\in P^*$.
We do not describe an antipode and a counit.
By this comultiplication, a tensor product of $\uq$-modules has
a $\uq$-module structure.

\nd
Let $*$ be the antiautomorphism of $\uq$ given by:
\beq
q^*=q,\q (q^h)^*=q^{-h},\q e_i^*=e_i,\q f_i^*=f_i.
\label{star}
\eeq
Let $\vee$ be the automorphism of $\uq$ given by:
\beq
q^{\vee}=q,\q (q^h)^{\vee}=q^{-h},
\q e_i^{\vee}=f_i,\q f_i^{\vee}=e_i.
\label{vee}
\eeq
These satisfy
\beq
**=\vee\vee={\rm id}, \qq *\vee=\vee*.
\label{star-vee}
\eeq

\subsection{Definition of Crystals}

Let us recall the definition of crystals \cite{K3,K4}.
The notion of a crystal
is motivated by abstracting the some combinatorial properties
of crystal bases.

\newtheorem{df2}{Definition}[section]
\begin{df2}
A {\it crystal} $B$ is a set with the following data:
\begin{eqnarray}
&&{\hbox{a map}}\q wt:B\lar P,\\
&&\vep_i:B\lar\ZZ\sqcup\{-\infty\},\q
  \vp_i:B\lar\ZZ\sqcup\{-\infty\},\q{\hbox{for}}\q i\in I,\\
&&\eit:B\lar B\sqcup\{0\},
\q\fit:B\lar B\sqcup\{0\}\q{\hbox{for}}\q i\in I.
\end{eqnarray}
Here 0 is an ideal element which is not included in $B$.
They are subject to the following axioms: For $b$,$b_1$,$b_2\in B$,
\begin{eqnarray}
&&\vp_i(b)=\vep_i(b)+\lan h_i,wt(b)\ran,\\
&&wt(\eit b)=wt(b)+\al_i{\hbox{ if }}\eit b\in B,\\
&&wt(\fit b)=wt(b)-\al_i{\hbox{ if }}\fit b\in B,\\
&&\eit b_2=b_1 {\hbox{ if and only if }} \fit b_1=b_2,
\label{eeff}\\
&&{\hbox{if }}\vep_i(b)=-\infty,
  {\hbox{ then }}\eit b=\fit b=0.
\end{eqnarray}
\end{df2}
{}From the axiom (\ref{eeff}),
we can consider the graph strucure on a crystal $B$.
\begin{df2}
\label{c-gra}
The crystal graph of crystal $B$ is
an oriented and colored graph given by
the rule : $b_1\mapright{i} b_2$ if and only if $b_2=\fit b_1$
$(b_1,b_2\in B)$.
\end{df2}
\begin{df2}
\label{df:mor}
\begin{enumerate}
\item
If $B$ has the weight decomposition $B=\bigsqcup_{\lm\in P}B_{\lm}$
where $B_{\lm}=\{b\in B|wt(b)=\lm\}$ for $\lm \in P$, we call
$B$ a $P$-weighted crystal.
\item
Let $B_1$ and $B_2$ be crystals.
A morphism of crystals $\psi:B_1\lar B_2$
is a map $\psi:B_1\sqcup\{0\}\lar B_2\sqcup\{0\}$
satisfying the following axioms:
\begin{eqnarray}
&&\hspace{-30pt}\psi(0)=0,
\label{psi(0)=0}\\
&&\hspace{-30pt}wt(b)=wt(\psi(b)),\q \vep_i(b)=\vep_i(\psi(b)),\q
\vp_i(b)=\vp_i(\psi(b))
\label{wt}\\
&&{\hbox{if }}b\in B_1{\hbox{ and }}\psi(b)\in B_2,\nonumber\\
&&\hspace{-30pt}\psi(\eit b)
=\eit\psi(b){\hbox{ if }}b\in B_1{\hbox{ satisfies }}
 \psi(b)\neq0{\hbox{ and }}\psi(\eit b)\neq0,\\
&&\hspace{-30pt}\psi(\fit b)
=\eit\psi(b){\hbox{ if }}b\in B_1{\hbox{ satisfies }}
 \psi(b)\neq0{\hbox{ and }}\psi(\fit b)\neq0.
\end{eqnarray}
\item
A morphism of crystals $\psi:B_1\lar B_2$
is called {\it strict} if the
associated map from $B_1\sqcup\{0\}\lar B_2\sqcup\{0\}$
commutes with all $\eit$ and $\fit$.
If $\psi$ is injective, surjective and strict,
$\psi$ is called an {\it isomorphism}.
\item
A crystal $B$ is a {\it normal},
if for any subset $J$ of $I$ such that
$((\al_i,\al_j))_{i,j\in J}$ is a positive symmetric matrix,
$B$ is isomorphic
to a crystal base of an integrable $U_q(\ge_{J})$-module, where
$U_q(\ge_{J})$ is the quantum algebra generated by $e_j$,
$f_j$ $(j\in J)$ and $q^h$ $(h\in P^*)$.
\end{enumerate}
\end{df2}

For crystals $B_1$ and $B_2$, we shall define their tensor product
$B_1\ot B_2$ as follows:
\begin{eqnarray}
&&B_1\ot B_2=\{b_1\ot b_2| b_1\in B_1 ,\, b_2\in B_2\},\\
&&wt(b_1\ot b_2)=wt(b_1)+wt(b_2),\\
&&\vep_i(b_1\ot b_2)={\hbox{max}}(\vep_i(b_1),
  \vep_i(b_2)-\lan h_i,wt(b_1)\ran),
\label{tensor-vep}\\
&&\vp_i(b_1\ot b_2)={\hbox{max}}(\vp_i(b_2),
  \vp_i(b_1)+\lan h_i,wt(b_2)\ran),
\label{tensor-vp}\\
&&\eit(b_1\ot b_2)=
\left\{
\begin{array}{ll}
\eit b_1\ot b_2 & {\mbox{ if }}\vp_i(b_1)\geq \vep_i(b_2)\\
b_1\ot\eit b_2  & {\mbox{ if }}\vp_i(b_1)< \vep_i(b_2),
\end{array}
\right.
\label{tensor-e}
\\
&&\fit(b_1\ot b_2)=
\left\{
\begin{array}{ll}
\fit b_1\ot b_2 & {\mbox{ if }}\vp_i(b_1)>\vep_i(b_2)\\
b_1\ot\fit b_2  & {\mbox{ if }}\vp_i(b_1)\leq \vep_i(b_2),
\label{tensor-f}
\end{array}
\right.
\end{eqnarray}

Here we understand that $0\ot b=b\ot 0=0$.
Let ${\cal C}(I,P)$ be the category
of crystals determined by the data $I$ and $P$.
Then $\ot$ is a functor
from ${\cal C}(I,P)\times{\cal C}(I,P)$
to ${\cal C}(I,P)$ and satisfies the
associative law:
$(B_1\ot B_2)\ot B_3\cong B_1\ot(B_2\ot B_3)$ by
$(b_1\ot b_2)\ot b_3\leftrightarrow b_1\ot (b_2\ot b_3)$.
 Therefore, the category of crystals is endowed
with the structure of tensor category.

For a crystal $B$, let $B^{\wedge}$ be the crystal given by
\beqn
&& B^{\wedge}:=\{b^{\wedge}|b\in B\},\nn\\
&& wt(b^{\wedge})=-wt(b),\qq \vep_i(b^{\wedge})=\vp_i(b),\qq
   \vp_i(b^{\wedge})=\vep_i(b),\nn\\
&& \eit(b^{\wedge})=(\fit b)^{\wedge},\qq
   \fit(b^{\wedge})=(\eit b)^{\wedge}.\nn
\eeqn
Then we have
\beq
(B_1\ot B_2)^{\wedge}\cong
B_2^{\wedge}\ot B_1^{\wedge}\qq {\rm by }\qq
(b_1\ot b_2)^{\wedge}\leftrightarrow b_2^{\wedge}\ot b_1^{\wedge}.
\label{b12}
\eeq

\newtheorem{ex}[df2]{Example}
\begin{ex}
\label{Example:crystal}
We give some examples of crystals.
\begin{enumerate}
\item
For $\lm\in P$, we set $T_{\lm}=\{t_{\lm}\}$ with
$$
wt(t_{\lm})=\lm,\q \vep_i(t_{\lm})=\vp_i(t_{\lm})=-\infty,\q
\eit(t_{\lm})=\fit(t_{\lm})=0.
$$
We can see that $T_{\lm}\ot T_{\mu}\cong T_{\lm+\mu}$ and
$B\ot T_{0}\cong T_{0}\ot B\cong B$ for any crystal $B$.
\item
For $i\in I$, we set $B_i:=\{(n)_i\,|\, n\in\ZZ\}$ with
\beqnn
&& wt((n)_i)=n\al_i,\qq \vep_i((n)_i)=-n,\qq \vp_i((n)_i)=n,\\
&& \vep_j((n)_i)=-\infty,\qq \vp_j((n)_j)
   =-\infty \q {\rm for }\q j\ne i,\\
&& \eit (n)_i=(n+1)_i,\qq \fit(n)_i=(n-1)_i,\\
&& \til e_j(n)_i=\til f_j(n)_i=0\q {\rm for }\q j\ne i.
\eeqnn
\item
For $\lm\in P_+$, let $(L(\lm),B(\lm))$ be the crystal base of a
$\uq$-integrable highest weight module $V(\lm)$.
$B(\lm)$ is the crystal associated with $V(\lm)$.
Set $B(-\lm)=B(\lm)^{\wedge}$. Then
$B(-\lm)$ is isomorphic to the crystal
associated with an integrable lowest
weight module $V(-\lm)$. For $b\in B(\lm)$,
$\vep_i(b)$ and $\vp_i(b)$ are given by
\beqnn
\vep_i(b)&=&
{\rm max}\{n|\eit^nb\ne 0\},\\
\vp_i(b)&=&
{\rm max}\{n|\fit^nb\ne 0\}.
\eeqnn
\item
Let $(L(\infty),B(\infty))$ be a crystal base of $\uqm$
$B(\infty)$ is a crystal associated with $\uqm$.
Set $B(-\infty)=B(\infty)^{\wedge}$.  $B(-\infty)$ is isomorphic to
the crystal associated with $U^+_q(\ge)$. In fact, we have
$B(-\infty)=B(\infty)^{\wedge}\cong B(\infty)^{\vee}$
by $b^{\wedge}\leftrightarrow b^{\vee}$.
For any $b\in B(\ify)$ and $i\in I$ there exists $k$ such that
$\eit^kb=0$. Then, $\vep_i(b)$ and $\vp_i(b)$ are given by
\beqnn
\vep_i(b)&=&
{\rm max}\{n|\eit^nb\ne 0\},\\
\vp_i(b)&=&\lan h_i,wt(b)\ran +\vep_i(b).
\eeqnn
\end{enumerate}
\end{ex}

%%%%%%%%% This section is Section 3 %%%%%%%%
\section{Crystals of modified quantum algebra}
\setcounter{equation}{0}
\renewcommand{\theequation}{\thesection.\arabic{equation}}
This section is devoted to
review \cite{K4},\cite{L2} (See also \cite{L1}).
\subsection{Modified quantum algebra and Crystal base}

For an integral weight $\lm\in P$,
let $\uq a_{\lm}$ be the left $\uq$-module
given by
\beq
\uq a_{\lm}:=\uq/\sum_{h\in P^*}\uq(q^h-q^{\lan h,\lm\ran}),
\label{ualm}
\eeq
where $a_{\lm}$ is the image of the unit by the canonical projection.
We set
\begin{equation}
\widetilde\uq=\bigoplus_{\lm\in P}\uq a_{\lm},
\label{eqn:def-util}
\end{equation}
which is called {\it modified quantum algebra}.

We shall see a crystal base of $\widetilde\uq$.
Taking $\lm\in P$ and
 choosing $\zeta,\mu\in P_+$ such that $\lm=\zeta-\mu$,
we get the following $\uq$-linear surjective homomorphism:
\beqn
\pi_{\zeta,\mu}  : \uq a_{\lm}&\longrightarrow &V(\zeta)\ot V(-\mu),
\label{pi}\\
               a_{\lm}&\mapsto &u_{\zeta}\ot u_{-\mu}.\nonumber
\eeqn
where $V(\zeta)$ and $V(-\mu)$ are
as in Example \ref{Example:crystal} (iii) and
$u_{\zeta}$ and $u_{-\mu}$ are their highest weight vector
and lowest weight vector respectively.

\newtheorem{thm3}{Theorem}[section]
\begin{thm3}[{\it cf} \cite{L2}]
\label{lusz}
For any $\lm\in P$,
there exists a unique pair $(L(\uq a_{\lm}),B(\uq a_{\lm}))$
which satisfies the following properties.
\begin{enumerate}
\item
We set $A:=\{f(q)\in\QQ(q)|{\hbox{$f$ has no pole at $q=0$}}\}$.
$L(\uq a_{\lm})$ is a free $A$-module
such that $\uq a_{\lm}\cong \QQ(q)\ot_A L(\uq a_{\lm})$ and
$B(\uq a_{\lm})$ is a $\QQ$-basis of
the $\QQ$-vector space $L(\uq a_{\lm})/qL(\uq a_{\lm})$.
\item
For any $\zeta,\mu\in P_+$ with $\lm=\zeta-\mu$, we have
$$
\pi_{\zeta,\mu} (L(\uq a_{\lm}))\subset L(\zeta)\ot_A L(-\mu),
$$
and the induced map $\bar\pi_{\zeta,\mu}$:
$$
\bar\pi_{\zeta,\mu}\,:\,L(\uq a_{\lm})/qL(\uq a_{\lm})\longrightarrow
(L(\zeta)/qL(\zeta))\ot(L(-\mu)/qL(-\mu)),
$$
satisfies $\bar\pi_{\zeta,\mu}
(B(\uq a_{\lm}))\subset B(\zeta)\ot B(-\mu)\,\sqcup \,\{0\}$.
\item
There is a structure of crystal
on $B(\uq a_{\lm})$ such that $\bar\pi_{\zeta,\mu}$ gives
a strict morphism of crystals
for any $\zeta,\mu\in P_+$ with $\lm=\zeta-\mu$.
\end{enumerate}
\end{thm3}

\vskip7pt
\nd
Set
$$
(L(\util),B(\util)):=
\bigoplus_{\lm\in P}(L(\uq a_{\lm}),B(\uq a_{\lm})).
$$
{\sl Remark.}
$B(\uq a_{\lm})$ is a normal crystal and
then $B(\util)$ is a normal crystal.

Let $B(\infty)$, $B(-\infty)$ and $T_{\lm}$ $(\lm\in P)$
be the crystals given in Example \ref{Example:crystal}.
The following theorem plays a significant role in this paper
(See \cite[Sec.3]{K4}).

\begin{thm3}
\label{U=BTB}
$B(\uq a_{\lm})\cong B(\infty)\ot T_{\lm}\ot B(-\infty)$ as a crystal.
\end{thm3}

\newtheorem{cor3}[thm3]{Corollary}
\begin{cor3}
$B(\wtil\uq)\cong \oplus_{\lm\in P}B(\infty)\ot T_{\lm}\ot B(-\infty)$
as a crystal.
\end{cor3}

\subsection{Description of the operation $*$}

By the definition of the operation $*$ given in (\ref{star}),
it acts on an element of $\util$ as follows;
Let $P$ be an element of $\uqm U_q^+(\ge)+U_q^+(\ge)\uqm$.
Arbitrary element $u\in U_q(\ge)a_{\lm}$ can be written in the form
$u=Pa_{\lm}$. Then we have
\beq
(Pa_{\lm})^*=a_{-\lm}P^*.
\label{star-on-u}
\eeq
Furthermore, we have the following results:
\begin{thm3}[\cite{K4}]
\label{thm-star-on-u}
\begin{enumerate}
\item
$L(\util)$ is invariant by $*$.
\item
$B(\util)^*=B(\util)$.
\item
For $\lm\in P$, $b_1\in B(\infty)$ and $b_2\in B(-\infty)$, we get
\beq
(b_1\ot t_{\lm}\ot b_2)^*=b_1^*\ot t_{-\lm-wt(b_1)-wt(b_2)}\ot b_2^*.
\label{description-star}
\eeq
\end{enumerate}
\end{thm3}

\subsection{Weyl group action and Extremal vectors}

This subsection is devoted to review \cite[Sec.7.8.9]{K4}.
Let $B$ be a normal crystal (See Definition \ref{df:mor} (iv)).
 Let us define the Weyl group action
on the underlying set $B$. For $i\in I$ and $b\in B$, we set
\begin{equation}
S_ib=
\left\{
\begin{array}{ll}
\fit^{\lan h_i,wt(b)\ran}b &{\rm if }\,\,
\lan h_i,wt(b)\ran\geq 0\\
\eit^{-\lan h_i,wt(b)\ran}b &{\rm if }\,\,
\lan h_i,wt(b)\ran<0.
\end{array}
\right.
\end{equation}

We can easily obtain the following formula:
\beq
 S_i^2={\hbox{id}},\qq
 S_i\eit=\fit S_i,\qq
 wt(S_i b)=s_i(wt(b)),
\label{33}
\eeq
where $s_i(\lm)=\lm-\lan h_i,\lm\ran\al_i$
is the simple reflection.

Let $\ge$ be a rank 2 finite dimensional Lie algebra,
and $W$ be the Weyl group associated with $\ge$.
We set $w_0=s_{i_1}\cdots s_{i_k}$ a reduced expression of
the longest element of $W$.
Here we get the following (\cite[Sec.7]{K4}):

\newtheorem{pro3}[thm3]{Proposition}
\begin{pro3}
Let $B$ be a normal crystal.
For any $b\in B$, $S_{i_1}\cdots S_{i_k}b$
does not depend on the choice of
reduced expression.
\end{pro3}

\begin{cor3}
$\{S_i\}_{i=1,2}$ satisfies the braid relation.
\end{cor3}

Thus for general $\ge$,
we know that $\{S_i\}_{i\in I}$ defines the Weyl group
action on a normal crystal $B$.
\newtheorem{df3}[thm3]{Definition}
\begin{df3}
\label{extremal-vector}
\begin{enumerate}
\item  Let $B$ be a normal crystal.
An element $b\in B$ is called $i$-extremal,
       if $\eit b=0$ or $\fit b=0$.
\item  An element $b\in B$ is called extremal if for any $l\geq0$,
$S_{i_1}\cdots S_{i_l}b$ is $i$-extremal for any $i$,
$i_1\cdots i_l\in I$.
\end{enumerate}
\end{df3}

\begin{thm3}
\label{conn-ext}
Any connected component of $B(\util)$ contains an extremal vector.
\end{thm3}

%%%%%%%%%%%% Section 4 %%%%%%%%%%%
\section{Criteria for Crystallized Peter-Weyl type decpomposition}
\setcounter{equation}{0}
\renewcommand{\theequation}{\thesection.\arabic{equation}}

\subsection{Right crystal structure on $\util$.}

Let $*$ be the antiautomorphism given in (\ref{star}).
By Theorem \ref{thm-star-on-u} (ii), we can define for $b\in B(\util)$;
\beqn
&& \vep_i^*(b):=\vep_i(b^*),\qq  \vp_i^*(b):=\vp_i(b^*),
\label{vep-star}\\
&& \eit^*(b):=(\eit (b^*))^*,\qq \fit^*(b):=(\fit (b^*))^*,
\label{ei-star}
\eeqn
By these, we can consider another crystal structure on $B(\util)$.
This another crystal structure has the following remarkable
property and that motivated us to consider
the Peter-Weyl type decomposition on $B(\util)$.

\newtheorem{thm4}{Theorem}[section]
\begin{thm4}[\cite{K4}]
\label{commute}
$\eit^*$ and $\fit^*$ are commutative with all $\eit$'s and $\fit$'s
on $B(\util)$, that is, on $B(\util)$ for any $i,j\in I$
$$
\til e_j\eit^*=\eit^*\til e_j,\q
\til f_j\eit^*=\eit^*\til f_j,\q
\til e_j\fit^*=\fit^*\til e_j,\q
\til f_j\fit^*=\fit^*\til f_j.
$$
\end{thm4}
A crystal endowed with another crystal structure as above
is called {\it bi-crystal}.

\subsection{Crystal $\bmax(\lm)$}
For $\lm\in P$, we set
\beq
\bmax(\lm)
:=\{b \in B(\uq a_{\lm})\,|\,{\hbox{$b^*$ is an extremal vector}}\}.
\label{bmax}
\eeq
By the following lemma, $\bmax(\lm)$ is a subcrystal of $B(\util)$.
\newtheorem{lem4}[thm4]{Lemma}
\begin{lem4}
\label{stability-bmax}
$$
\eit \bmax(\lm)\subset \bmax(\lm)\sqcup\{0\},\qq
\fit \bmax(\lm)\subset \bmax(\lm)\sqcup\{0\}.
$$
\end{lem4}

\nd
{\sl Proof.\,}
Let $b$ be an element of $\bmax(\lm)$.
For any $i,i_1,\cd, i_k\in I$, by the definition of extremal vector,
$b^*$ satisfies,
\beq
\eit S_{i_1}\cd S_{i_k}b^*=0 \q {\hbox{ or }}\q
\fit S_{i_1}\cd S_{i_k}b^*=0.
\label{e-star}
\eeq
By operating $*$ on the both sides of (\ref{e-star})
and by the fact $**={\rm id}$,
it follows that
\beq
\eit^* S_{i_1}^*\cd S_{i_k}^*b=0 \q {\hbox{ or }}\q
\fit^* S_{i_1}^*\cd S_{i_k}^*b=0,
\label{stra-e-star}
\eeq
where $S_i^*b:=(S_i(b^*))^*$.
This (\ref{stra-e-star}) is written
by using only $\eit^*$'s and $\fit^*$'s.
Thus, by Theorem \ref{commute}, for any $j\in I$ we obtain
$$
\eit^* S_{i_1}^*\cd S_{i_k}^*\til e_jb=0 \q {\hbox{ or }}\q
\fit^* S_{i_1}^*\cd S_{i_k}^*\til e_jb=0
$$
and
$$
\eit^* S_{i_1}^*\cd S_{i_k}^*\til f_jb=0 \q {\hbox{ or }}\q
\fit^* S_{i_1}^*\cd S_{i_k}^*\til f_jb=0.
$$
Therefore, by operating $*$ on the both sides, we get
$$
\eit S_{i_1}\cd S_{i_k}(\til e_jb)^*=0 \q {\hbox{ or }}\q
\fit S_{i_1}\cd S_{i_k}(\til e_jb)^*=0
$$
and
$$
\eit S_{i_1}\cd S_{i_k}(\til f_jb)^*=0 \q {\hbox{ or }}\q
\fit S_{i_1}\cd S_{i_k}(\til f_jb)^*=0.
$$
We have completed the proof.\qed

\nd
The crystal $\bmax(\lm)$ has the following properties:
\newtheorem{pr4}[thm4]{Proposition}
\begin{pr4}[\cite{K4}]
\label{pro-bmax}
\begin{enumerate}
\item
For any $i\in I$, $S_i^*$ gives an isomorphism;
\beqn
S_i^*: \bmax(\lm) & \mapright{\sim} & \bmax(s_i\lm)
\label{si-iso}\\
b &\longmapsto & S_i^*b\nonumber
\eeqn
\item
Let $u_{\lm}\in B(\uq a_{\lm})$ be the corresponding element to
$u_{\infty}\ot t_{\lm}\ot u_{-\infty}$ through the isomorphism
$B(\uq a_{\lm})\cong B(\infty)\ot T_{\lm}\ot B(-\infty)$.
Then $u_{\lm}\in \bmax(\lm)$.
\end{enumerate}
\end{pr4}

\subsection{Criteria for Peter-Weyl type decomposition}

For $\lm\in P$, set
\beqnn
B(\lm) & := &
\{\til X_{i_1}\cd \til X_{i_k}u_{\lm}\,;\, \til X_i=\til e_i,
\til f_i, \,\,i_j\in I\}
\setminus\{0\}
\subset B(\uq a_{\lm}),\\
B(\lm)^* & := &
\{\til X_{i_1}^*\cd \til X_{i_k}^*u_{\lm}^*\,;\, \til X_i
=\til e_i,\til f_i, \,\,i_j\in I\}
\setminus\{0\}.
\eeqnn
{\sl Remark.} \,
\begin{enumerate}
\item
In terms of crystal graph, $B(\lm)$ is a connected
component of $B(\uq a_{\lm})$ containing $u_{\lm}$.
\item
If $\lm$ is a dominant, $B(\lm)$ coincides with
the one in Example \ref{Example:crystal} (iii).
\item
$B(\lm)^*$ is stable by the actions of $\eit^*$ and $\fit^*$, that is,
\beq
\eit^* B(\lm)^*\subset B(\lm)^*\sqcup\{0\},\qq
\fit^* B(\lm)^*\subset B(\lm)^*\sqcup\{0\}.
\label{stab-bstar}
\eeq
\item
By the definition of $\eit^*$ and $\fit^*$, we have
$$
\til X_{i_1}^*\cd \til X_{i_k}^*u_{\lm}^*
=(\til X_{i_1}\cd \til X_{i_k}u_{\lm})^*.
$$
Therefore, we get
\beq
B(\lm)^*=\{b^*\,;\,b\in B(\lm)\}.
\label{bstar=bstar}
\eeq
\end{enumerate}

We consider the following three conditions.
\begin{description}
\item[(C1)]
For any extremal vector $b\in B(\util)$, there exists
an embeding of crystal
$$
B(wt(b))\hookrightarrow B(\util),
$$
given by $u_{wt(b)}\mapsto b$.
\item[(C2)]
For any $\lm\in P$, $B(\lm)_{\lm}=\{u_{\lm}\}$.
\item[(C3)] (transitivity of extremal vectors)\,\,
For any extremal vectors $b_1,\,b_2\in B(\lm)$,
there exist $i_1,\cd i_k$ such that
$$
b_2=S_{i_1}\cd S_{i_k}b_1.
$$
\end{description}
\begin{pr4}
\label{Peter-Weyl}
If $\util$ satisfies the conditions (C1), (C2) and (C3),
there exists the following isomorphism of bi-crystal;
\beq
B(\util)\cong
\bigoplus_{\lm\in P/W}\bmax(\lm)\ot B(-\lm)^*,
\label{P-W-formula}
\eeq
where $W$ is the Weyl group associated with $\ge$ and
an isomorphism of bi-crystal is, by definition, an isomorphism for
both crystal structures.
\end{pr4}
{\sl Remark.\,}
The tensor product in (\ref{P-W-formula}) has a different meaning
from usual tensor product of crystals.
In the tensor product in the R.H.S. of (\ref{P-W-formula}),
$\eit$ and $\fit$ act on the left component
and $\eit^*$ and $\fit^*$ act on the right component, that is,
for $u\ot v\in \bmax(\lm)\ot B(-\lm)^*$,
$\til X_i(u\ot v)=\til X_iu\ot v$  and
$\til X_i^*(u\ot v)=u\ot \til X_i^*v$ ($X_i=e_i,f_i$).

\vskip5pt
\nd
{\sl Proof.\,}
In order to show the proposition we shall see the following lemmas.
\begin{lem4}
\label{lem-iso}
We set
\beq
\wtil B(\lm):=
\{\til X_{i_1}^*\cd \til X_{i_k}^*b\,;\,b\in\bmax(\lm),\,\,
i_j\in I,\,\,X_i=e_i,f_i\}\setminus\{0\}.
\label{btil}
\eeq
If the conditions (C1) and (C2) hold,
$\wtil B(\lm)$ has a bi-crystal structure and we have the
following isomorphism of bi-crystal;
\beqn
\vp:\bmax(\lm)\ot B(-\lm)^* & \mapright{\sim} &  \wtil B(\lm),
\label{small-iso}\\
b\ot \til X_{i_1}^*\cd \til X_{i_k}^*u_{-\lm}^*
&\mapsto & \til X_{i_1}^*\cd \til X_{i_k}^*b.
\eeqn
\end{lem4}
\begin{lem4}
\label{blm=blm'}
Assume that the conditions (C3) holds.
The following (A) and (B) are equivalent:
\begin{description}
\item[(A)]
$\lm,\lm'\in P$ satisfy $\lm'=w\lm$ for some $w\in W$,
\item[(B)]
$\wtil B(\lm)=\wtil B(\lm').$
\end{description}
\end{lem4}

\begin{lem4}
\label{emp}
We assume that (C3) holds.
For $\lm,\lm'\in P$ if there is no $w\in W$ such that
$\lm'=w\lm$, we have
$$
\wtil B(\lm)\cap \wtil B(\lm')=\emptyset.
$$
\end{lem4}

\vskip5pt
\nd
{\sl Proof of Lemma \ref{lem-iso}.\,}
To see that $\wtil B(\lm)$ is a well-defined bi-crystal,
it is sufficient to show the
stability of $\wtil B(\lm)$ by the actions of $\eit$, $\fit$,
$\eit^*$ and $\fit^*$, which is derived easily from
Theorem \ref{commute}
and Lemma \ref{stability-bmax}.

Next we shall see the well-definedness of the map $\vp$.
Since the definition of the map $\vp$ depends on the expression of
a right component, we shall show
\begin{description}
\item[(I)]
If $X_{i_1}^*\cd X_{i_k}^*u_{-\lm}^*=0$,
$X_{i_1}^*\cd X_{i_k}^*b=0$ for any $b\in \bmax(\lm)$.
\item[(II)]
If $X_{i_1}^*\cd X_{i_k}^*u_{-\lm}^*
=X_{j_1}^*\cd X_{j_l}^*u_{-\lm}^*\ne 0$,
 then
$$
X_{i_1}^*\cd X_{i_k}^*b
=X_{j_1}^*\cd X_{j_l}^*b\ne 0
$$
for any $b\in \bmax(\lm)$.
\end{description}

\nd
(I)\,\,\,Under the assumption
$(X_{i_1}^*\cd X_{i_k}^*u_{-\lm}^*)^*=
X_{i_1}\cd X_{i_k}u_{-\lm}=0$
it is enough to show
$$
(X_{i_1}^*\cd X_{i_k}^*b)^*=X_{i_1}\cd X_{i_k}b^*=0.
$$
By the definition of $\bmax(\lm)$ and
Theorem \ref{thm-star-on-u} (iii),
$b^*$ is an extremal vector and has a weight $-\lm$.
Let $B'$ be the connected component containing $b^*$.
Since the condition (C1) holds,
$B'\cong B(-\lm)$ by $b^*\leftrightarrow u_{-\lm}$.
Thus, we obtain (I).

\nd
(II)\,\,\,
We set $\til Y_i=\eit$ if $\til X_i=\fit$
and $\til Y_i=\fit$ if $\til X_i=\eit$.
Then we have
$\til Y_{i_k}^*\cd \til Y_{i_1}^*\til X_{i_1}^*\cd \til X_{i_k}^*b=b$.
Therefore, it is sufficent to show that
\beq
\til Y_{i_k}^*\cd \til Y_{i_1}^*\til X_{j_1}^*\cd \til X_{j_l}^*b=b.
\label{b=b}
\eeq
We know that $b^*$ is an extremal vector and $wt(b^*)=-\lm$.
Let us denote $B'$ for the connected component
containing $b^*$.
By applying $*$ on the L.H.S of (\ref{b=b}), we get
$$
(\til Y_{i_k}^*\cd \til Y_{i_1}^*\til X_{j_1}^*\cd \til X_{j_l}^*b)^*
=\til Y_{i_k}\cd \til Y_{i_1}\til X_{j_1}\cd \til X_{j_l}b^*\in B'
$$
Since $wt(X_{i_1}^*\cd X_{i_k}^*u_{-\lm}^*)
=wt(X_{j_1}^*\cd X_{j_l}^*u_{-\lm}^*)$, we get

\beq
wt(\til Y_{i_k}\cd \til Y_{i_1}\til X_{j_1}\cd \til X_{j_l}b^*)
=wt(b^*)=-\lm.
\label{wt}
\eeq
By the condition (C2) and the fact $B'\cong B(-\lm)$,
we know that $B'_{-\lm}=\{b^*\}$.
Thus, by (\ref{wt}) we have
\beq
\til Y_{i_k}\cd \til Y_{i_1}\til X_{j_1}\cd \til X_{j_l}b^*
=b^*.
\label{YXb=b}
\eeq
Applying $*$ on the both sides of (\ref{YXb=b}),
we obtain (\ref{b=b}).

By Theorem \ref{commute}, it is trivial that
$$
\vp\circ\til X_i=\til X_i\circ\vp \qq{\rm and }\qq
\vp\circ\til X_i^*=\til X_i^*\circ\vp
$$
Thus, the map $\vp$ is a morphism of bi-crystal.

We shall see the surjectivity of the map $\vp$.
It is enough to show that
if $\til X_{i_1}^*\cd \til X_{i_k}^*b \in \wtil B(\lm)$
for $b\in \bmax(\lm)$,
$\til X_{i_1}^*\cd \til X_{i_k}^*u_{-\lm}^*\ne 0$.
It is shown by the fact that $B'\cong B(-\lm)$,
where $B'$ is the connected component including $b^*$.
%%%%%%%%
%By the assumption, we have
%\beq
%0\ne (\til X_{i_1}^*\cd \til X_{i_k}^*b)^*
%=\til X_{i_1}\cd \til X_{i_k}b^*.
%\label{41}
%\eeq
%Since the vector $b^*$ is an extremal vector and $wt(b^*)=-\lm$,
%by the condition (C1) and (\ref{41}),
%we get $\til X_{i_1}\cd \til X_{i_k}u_{-\lm}\ne 0$.
%Then applying $*$, we obatin
%$$
%0\ne (\til X_{i_1}\cd \til X_{i_k}u_{-\lm})^*
%= \til X_{i_1}^*\cd \til X_{i_k}^*u_{-\lm}^*.
%$$

Finally, we shall see the injectivity of the map $\vp$.
For $u=b_1\ot \til X_{i_1}^*\cd \til X_{i_k}^*u_{-\lm}^*$ and
$v=b_2\ot \til X_{j_1}^*\cd \til X_{j_l}^*u_{-\lm}^*\in
\bmax(\lm)\ot B(-\lm)^*$,
we shall prove that
if $\vp(u)=\vp(v)$, $u=v$.
Applying $*$ on $\vp(u)=\vp(v)$, we get
\beq
\til X_{i_1}\cd \til X_{i_k}b_1^*=
\til X_{j_1}\cd \til X_{j_l}b_2^*.
\label{42}
\eeq
Since both
$b_1^*$ and $b_2^*$ are extremal vectors with
weight $-\lm$ and
(\ref{42}) implies that $b_1^*$ and $b_2^*$ are contained in the
same connected component,
by the conditions (C1) and (C2) we obtain $b_1^*=b_2^*$ and
$\til X_{i_1}\cd \til X_{i_k}u_{-\lm}^*=
\til X_{j_1}\cd \til X_{j_l}u_{-\lm}^*$.
Thus we get $u=v$.
Now we have completed the proof of Lemma \ref{lem-iso}.
\qed

\vskip6pt
\nd
{\sl Proof of Lemma \ref{blm=blm'}.\,}
In order to see (A) $\Rightarrow$ (B),
it is enough
to show that for $\lm'=s_i\lm$
\beq
\wtil B(\lm')\subset \wtil B(\lm).
\label{43}
\eeq
By Proposition \ref{pro-bmax} (i),
$S_i^*:\bmax(\lm)\mapright{\sim} \bmax(\lm')$.
Then we have
$$
\bmax(\lm')=S_i^*\bmax(\lm)\subset \wtil B(\lm).
$$
Therefore, we get
$$
\wtil B(\lm')\subset \wtil B(\lm).
$$
Since $\lm=s_i\lm'$, we also get $\wtil B(\lm)\subset \wtil B(\lm')$.

Next we shall see that (B) $\Rightarrow$ (A).
For $u_{\lm}\in \bmax(\lm)$ there exists $b\in \bmax(\lm')$
and ${i_1},\cd ,{i_k}$
such that
$$
u_{\lm}=\til X_{i_1}^*\cd \til X_{i_k}^*b.
$$
Applying $*$ on the both sides, we get
$$
u_{\lm}^*=\til X_{i_1}\cd \til X_{i_k}b^*.
$$
This implies that $u_{\lm}^*$ and $b^*$ belong
to the same connected component.
Since both $u_{\lm}^*$ and $b^*$ are
extremal vectors,  by (C3) there exist $j_1,\cd,j_l$
such that
\beq
u_{\lm}^*=S_{j_1}\cd S_{j_l}b^*.
\label{44}
\eeq
Here note that $wt(u_{\lm}^*)=-\lm$ and
$wt(b^*)=-\lm'$ by Theorem \ref{thm-star-on-u} (iii).
Thus, by (\ref{33}) and
(\ref{44}) we obtain $\lm=s_{j_1}\cd s_{j_l}\lm'$.
\qed

\vskip5pt
\nd
{\sl Proof of Lemma \ref{emp}}\,\,
We assume that $\wtil B(\lm)\cap \wtil B(\lm')\ne\emptyset$ and
$b\in\wtil B(\lm)\cap \wtil B(\lm')$. There exist
$i_1,\cd,i_k,j_1,\cd,j_l\in I$ such that
\beqnn
&&b_1:=\til X_{i_1}^*\cd \til X_{i_k}^*b\in \bmax(\lm),\\
&&b_2:=\til X_{j_1}^*\cd \til X_{j_l}^*b\in \bmax(\lm').
\eeqnn
By the definition of $\bmax(\lm)$,
$b_1^*$ and $ b_2^*$ are extremal vectors and
\beq
wt(b_1^*)=-\lm \qq{\rm and}\qq wt(b_2^*)=-\lm'.
\label{20}
\eeq
We also get
\beq
b_2=\til X_{j_1}^*\cd\til X_{j_l}^*\til Y_{i_k}^*\cd\til Y_{i_1}^*b_1,
\label{21}
\eeq
where $\til Y_i$ is the same one as
in the proof of Lemma \ref{lem-iso}.
Applying * on the both sides of (\ref{21}),
we get
$$
b_2^*=\til X_{j_1}\cd\til X_{j_l}\til Y_{i_k}\cd\til Y_{i_1}b_1^*.
$$
This implies that
$b_1^*$ and $b_2^*$ are in the same connected component.
By virtue of (C3),
there exist $a_1,\cd, a_n$ such that
\beq
b_2^*=S_{a_1}\cd S_{a_n} b_1^*.
\label{211}
\eeq
By (\ref{33}), (\ref{20}) and (\ref{211}), we obtain
$\lm'=s_{a_1}\cd s_{a_n}\lm$, which is a contradiction.\qed

\vskip7pt
\nd
{\sl Proof of Proposition \ref{Peter-Weyl}.\,}
By Lemma \ref{blm=blm'}, we know that
for $\lm,\lm'\in P$ if there exists $w\in W$ such that
$\lm'=w\lm$, $\wtil B(\lm)=\wtil B(\lm')$ and
otherwise, $\wtil B(\lm)\cap \wtil B(\lm')=\emptyset$
by Lemma \ref{emp}.
Therefore, we get
$\sum_{\lm\in P}\wtil B(\lm)=\oplus_{\lm\in P/W}\wtil B(\lm)$
and then
\beq
\bigoplus_{\lm\in P/W}\wtil B(\lm)
\subset
B(\util).
\label{hook}
\eeq
On the other hand,
for any $b\in B(\util)$ let $B'$ be the
connected component containing $b^*$.
By Theorem \ref{conn-ext},
there exist $\til X_{i_1},\cd,\til X_{i_k}$
such that $\til X_{i_1}\cd\til X_{i_k}b^*$
is an extremal vector. Then there exists some
$\mu\in P$ such that
$$
(\til X_{i_1}\cd\til X_{i_k}b^*)^*
=\til X_{i_1}^*\cd\til X_{i_k}^*b\in \bmax(\mu).
$$
This implies $b\in \widetilde B(\mu)$.
Therefore, we get
\beq
B(\util)\subset \sum_{\lm\in P}\wtil B(\lm)=
\bigoplus_{\lm\in  P/W}
\wtil B(\lm).
\label{subset}
\eeq
By (\ref{hook}) and (\ref{subset}) we obtain
$$
B(\util)=\bigoplus_{\lm\in  P/W}
\wtil B(\lm).
$$
Finally, by
Lemma \ref{lem-iso}, we get the desired result.    \qed

%%%%%%%%%%%%%%%  Section 5  %%%%%%%%%%%%%%
\section{Crystallized Peter-Weyl type
decomposition for level 0 part of $B(\til U_q(\slh))$}
\setcounter{equation}{0}
\renewcommand{\theequation}{\thesection.\arabic{equation}}

In this secton we set $\ge=\slh$.
As an application of Proposition \ref{Peter-Weyl},
we shall give the Peter-Weyl type decomposition for
the level 0 part of $B(\wtil U_q(\slh))$ explicitly.
The notations and terminologies using in this section
follow \cite{K3},  \cite[Sec5--Sec7]{N}.

\subsection{Crystals of level 0 part of $\wtil U_q(\slh)$}
In this subsection we shall review \cite{K3},\cite{N}.

\vskip5pt
Let us denote $\util_0$ for the level 0 part of $\util$ given by
$$
\util_0:=\{b\in\util|\lan c,wt(b)\ran=0\}
=\bigoplus_{\lm\in P_0:=\{\lm\in P\,|\,\lan c,\lm\ran=0\}}\uq a_{\lm}.
$$

\nd
We set
$$
B_{\ify}:=\{(n)|n\in\ZZ\},
\qq (wt(n)=2n(\Lm_0-\Lm_1)).
$$
We define the crystal structure on $B_{\ify}$ by
\beqnn
&&\til e_1(n)=\til f_0(n)=(n-1),\qq\qq
\til e_0(n)=\til f_1(n)=(n+1),\\
&&\vep_1(n)=\vp_0(n)=n,\qq\qq
\vep_0(n)=\vp_1(n)=-n.
\eeqnn
{\sl Remark.}
We shall identify $B_{\ify}$ with $\ZZ$. Then we can consider
summation, subtraction and absolute value for elements in $B_{\ify}$.

 We set
\begin{eqnarray}
  &&{\cal P}(\infty):=\{(\cdot\cdot,i_k,i_{k+1},\cdot\cdot,i_{-1})|
  i_k\in B_{\infty}{\hbox{ and if }} |k|\gg0, i_k=(0)\},\\
  &&{\cal P}(-\infty):=\{(i_0,\cdot\cdot,i_k,i_{k+1},\cdot\cdot)|
  i_k\in B_{\infty}{\hbox{ and if }} |k|\gg0, i_k=(0)\},
\end{eqnarray}

Then we get the following isomorphisms of crystal
(\cite{KKM},\cite[5.1]{N})
$$
B(\ify)\cong {\cal P}(\ify),\qq
B(-\ify)\cong {\cal P}(-\ify).
$$
For an integer $m$ we set
\beq
\begin{array}{rcl}
&& {\cal P}_m :=
\left\{
p=(\cdots,i_k,i_{k+1},\cdots,i_{-1},i_0,i_1,\cdots,i_l,i_{l+1},\cdots)
|\right. i_k\in B_{\infty}  \\
&&\,\,\,\left.{\hbox{ if }}k\ll0,\,\,i_k=(0)
  {\hbox{ and if }}l\gg0,\,\, i_{2l}
=(m){\hbox{ and }}i_{2l+1}=(-m)\right\}.\nonumber
\end{array}
\label{path}
\eeq
Furthermore, let us associate a weight with $p\in{\cal P}_m$
by the following:
\beq
\begin{array}{rcl}
\lefteqn{wt(p)=(\sum_{k\in\ZZ}i_{k-1}+i_{k})(\Lm_0-\Lm_1) } \\
&& \qq  +(l+\sum_{k\in\ZZ}k({\rm max}\{i_{k-1},-i_{k}\}
-{\rm max}\{g_{k-1}, -g_k\}))\del,
 \end{array}
  \label{path-wt}
 \eeq
where $l$ is an integer and
$(g_k)_k$ is the element in ${\cal P}_m$ given by
$g_k=(0)$ $(k<0)$ and $g_k=((-)^km)$ $(k\geq 0)$.
We set ${\cal P}_{m,l}={\cal P}_m$ as a set
and weight of an element in ${\cal P}_{m,l}$
 is given by (\ref{path-wt}).
For $\lm=m(\Lm_0-\Lm_1)+l\del$ we have
\beq
{\cal P}_{m,l}\cong B(U_q(\slh)a_{\lm}).
\label{544}
\eeq
Let us call an element of ${\cal P}_{m,l}$ a {\it path}.

For $n\in \ZZ_{\geq 0}$, let ${\cal P}_m(n)$ be the subset of
${\cal P}_m$ with $n$ walls
as in \cite[Sec 5.3]{N}.
For a sequence in $m$-domain configuration (see \cite[Sec 7.1]{N})
$\vec t=(t_1,\cd,t_{n-1})$ and a sequence
$\vec c=(c_1,\cd,c_{n-1})\in \ZZ^{n-1}_{\geq0}$
let $\mlpathtc$ be the same object as in \cite[Sec 7.4]{N}.
By \cite{N} we get the following
\newtheorem{thm5}{Theorem}[section]
\begin{thm5}
$\mlpathtc$ is a connected component of $\mlpath$ and
any connected component of $\mlpath$ coincides with
some $\mlpathtc$.
\end{thm5}

\nd
The crystal $B(\ify)$ has another path realization. (See \cite{K3}).
For $i\in I$, let $B_i$ be the crystal as in
Example \ref{Example:crystal} (ii).
We define the map
$\Psi_i:B(\infty)\longrightarrow B(\infty)\ot B_i$
by
$\Psi_i(b)=b_0\ot \fit^m(0)_i$,
where $m=\vep_i^*(b)$ and $b_0=\eit^{*m}b$.
Then we have the following;
\begin{thm5}[\cite{K3}]
\label{B=BBi}
For any $i\in I$, $\Psi_i$ gives a strict embeding of crystals.
\end{thm5}
Let us take a sequence $i_1,i_2,\cd,i_n\in I=\{0,1\}$ such that
$i_k\ne i_{k+1}$.

\nd
{\sl Remark.\,}
In the case $\ge=\slh$, there are only two choices of
a sequence $i_1,i_2,\cd$,
$(i_1,i_2,\cd)=(1,0,1,0,\cd)$ or $(0,1,0,1,\cd)$.

\nd
 By iterating $\Psi_i$ we obtain;
\beqn
&&\Psi_{i_n,\cd,i_1}:B(\ify)\mapright{\Psi_{i_1}}
B(\ify)\ot B_{i_1}\mapright{\Psi_{i_2}\ot{\rm id}}
B(\ify)\ot B_{i_2}\ot B_{i_1}\\
&&\qq \rightarrow\cd\rightarrow
B(\ify)\ot B_{i_n}\ot\cd B_{i_1}.
\label{bi-path}
\eeqn
For any $b\in B(\ify)$ if we take $n\gg0$, $\Psi_{i_n,\cd,i_1}(b)$
can be written in the following form;
$$
\Psi_{i_n,\cd,i_1}(b)=u_{\ify}\ot \til f_{i_n}^{a_n}(0)_{i_n}\ot
\cd\ot\til f_{i_1}^{a_1}(0)_{i_1}\qq (a_k\geq 0).
$$
For $m>n$, since $\Psi_{i_m,\cd,i_1}(b)=
u_{\ify}\ot(0)_{i_m}\ot \cd\ot (0)_{i_{n+1}}\ot
\til f_{i_n}^{a_n}(0)_{i_n}\ot
\cd\ot\til f_{i_1}^{a_1}(0)_{i_1},$ the sequence
$a_1,a_2,\cd$ does not depend on the choice of $n$.
Thus, we have the follwoing embeding of crystal;
\beq
\Psi:B(\ify)\longrightarrow
\{(\cd,a_n,a_{n-1},\cd,a_1)|\,\,a_n\in\ZZ_{\geq 0}, \q a_k=0\q(k\gg 0)\},
\label{psi}
\eeq
where $a_k$ means ${\til f_{i_k}}^{a_k}(0)_{i_k}$.
In \cite{K3}, the image of the map $\Psi$ is described explicitly
for any rank 2 Kac-Moody Lie algebra. In paticular, for $\ge=\slh$
we have
\newtheorem{pr5}[thm5]{Proposition}
\begin{pr5}
\label{b-inf-sl2}
We set $z_n={n\over{n-1}}$ ($n\geq 2$). Then we have
\beq
{\rm Im}\,\Psi=
\left\{(\cd,a_n,a_{n-1},\cd,a_1)|
\begin{array}{ll}
&a_n\in\ZZ_{\geq 0} \q{\rm and  }\q a_k=0 \q(k\gg 0)\\
&a_{n+1}\leq z_na_n\q(n\geq 2)
\end{array}
\right\}.
\eeq
\end{pr5}

\subsection{Description of $\bmax(\lm)$ with level$(\lm)=0$}

For $b=b_1\ot t_{\lm}\ot b_2\in B(\ify)\ot T_{\lm}\ot B(-\ify)$,
by Theorem \ref{thm-star-on-u} (iii) we get
$$
b^*=b_1^*\ot t_{-\lm-wt(b_1)-wt(b_2)}\ot b_2^*
=b_1^*\ot t_{-wt(b)}\ot b_2^*\in B(\ify)\ot T_{-wt(b)}\ot B(-\ify).
$$
Thus, by (\ref{bmax}) and Theorem \ref{U=BTB} we can write;
\beq
\bmax(\lm)=
\{b^*|b\in B(\util)_{-\lm}{\hbox{ and $b$ is an extremal vector}}\}
\label{bmax2}
\eeq
For level 0 weight $\lm$, we shall describe $\bmax(\lm)$ explicitly.
Set $\lm=m(\Lm_0-\Lm_1)+l\del$
where $l,m\in \ZZ$.
Let $p\in {\cal P}_k$ be an extremal vector.
Since  all walls (see \cite[Sec 5.3]{N})
in $p$ are simultaneously $+$ or $-$
by \cite[Theorem 7.18]{N}, we have
\beq
wt(p)=
\left\{
\begin{array}{ll}
n(\Lm_0-\Lm_1)+j\del&{\rm \q if \,\, all\,\, walls\,\, are\q}+,\\
n(\Lm_1-\Lm_0)+j\del&{\rm \q if\,\, all\,\, walls\,\, are\q}-,
\end{array}
\right.
\label{m01}
\eeq
where $n$ is a number of walls in $p$ and $j$ is some integer.
Thus, we obtain the following lemma.
\newtheorem{lem5}[thm5]{Lemma}
\begin{lem5}
\label{mmm}
We set ${\lm}=m(\Lm_0-\Lm_1)+l\del$ with $m\geq 0$
$($resp. $m\leq 0$$)$.
If  $b=(\cd,i_{k},i_{k+1},\cd)$ is an extremal vector
with the weight $-\lm$.
then $b$ has $|m|$ walls and all walls
in $b$ are $-$ $($resp. $+$$)$
and $b^*\in {\cal P}_{m,l}$,
{}.
\end{lem5}

We shall describe the operation $*$ on an extremal vector
$p=(\cd,i_{k},i_{k+1},\cd)\in B(\util_0)$ more precisely.
Set $b=(\cd,i_{k},i_{k+1},\cd,i_{-1})
\in {\cal P}(\ify)\cong B(\ify)$.
We can define walls and domains
in $b$ by the similar way to the one in
\cite[Sec 5.3]{N}
\newtheorem{def5}[thm5]{Definition}
\begin{def5}
\begin{enumerate}
\item
For $b=(\cd,i_k,i_{k+1},\cd,i_{-1})\in B(\ify)$, by definition,
there are $+$ $($resp. $-$$)$
walls at position $k\leq -1$
if $i_{k-1}+i_k>0$ $($resp. $i_{k-1}+i_k<0$$)$
and $|i_{k-1}+i_k|$ is called the number
of walls at position $k$. We also call
$\sum_{k\leq -1}|i_{k-1}+i_k|$ the total
number of walls in $b$.
\item
A {\sl domain} in $b$ is a finite subsequence
$i_a,i_{a+1},\cd,i_b$ such that $i_{a-1}+i_a\ne 0$, $i_b+i_{b+1}\ne 0$
and $i_j+i_{j+1}=0$ for $a\leq j<b$, namely, a subsequence surrounded
by neighboring two walls
or a subsequence $\cd,i_{c-1},i_c$ $($resp. $i_{c},i_{c+1},\cd$$)$
such that $i_c+i_{c+1}\ne 0$ and $i_{j-1}+i_j=0$ for $j\leq c$
$($resp. $i_{c-1}+i_c\ne 0$ and $i_{j-1}+i_j=0$ for $j>c$$)$, namely,
a subsequence of $p$ on the left $($resp. right$)$ of
the left-most $($resp. right-most$)$ wall or
an empty subsequence without entry which appears
between two neighboring walls in the same position.

\end{enumerate}
\end{def5}
{\sl Remark.}\,\,
The left-most domain is an infinite domain but the right-most domain
is not infinite.

\newtheorem{ex5}[thm5]{Example}
\begin{ex5}
For a path $p=(\cd,0,0,1,-1,3,-3)$,
there is a  zero-length domain between
entries $-1$ and $3$ and there are one infinite domain
$\cd,0,0,0,0$ and two finite domains $1, -1$ and $3, -3$.
\end{ex5}
%The left-most domain is the infinite
%subsequence in $b$ given by
%$$
%\cd,i_{m-2},i_{m-1}
%$$
%where $m$ is the smallest integer such that $i_m\ne0$.
%Thus, any entry in the left most domain is 0.
%The right-most domain is the subsequence
%$i_k,\cd,i_{-1}$
%where $k\leq -1$ is the lagrgest integer
%such that $i_{k-1}+i_k\ne 0$.

\vskip5pt
\nd
If the total number of walls is $n$,
there are $n+1$ domains in $b$ and
we shall denote $d_0, d_1,\cd, d_n$ for them, where
$d_0$ is the left-most domain and $d_n$ is the right-most domain.
All domains except $d_0$ include finite number of entries.
For a domain $d_j$ $(j>0)$, let $l_j$ be the length of
$d_j(=$ number of entries in $d_j)$.
In particular, the total sum of lenghts
$$
l(b):=\sum_{j=1}^{n}l_j
$$
is called {\it length of } $b$.

\begin{pr5}
\label{b1}
Let $p=(\cd,i_k,i_{k+1},\cd)\in B(\util_0)$ be an extremal vector
and set $b=(\cd,i_k,\cd,i_{-2},i_{-1})\in {\cal P}(\ify)$.
Suppose that $b$ has $n$ walls and
set $l'_j=\sum_{i=1}^jl_i$ $(l'_n=l(b))$. Then we have
\beq
b=
\underbrace{\til f_{j_{l(b)}}^{n}
\cd\til f_{j_{l'_{n-1}+1}}^{n}}_{l_{n}}
\cd\cd
\underbrace{\til f_{j_{l'_2}}^{2}
\cd\til f_{j_{l'_1+1}}^{2}}_{l_{2}}
\underbrace{\til f_{j_{l'_1}}
\cd\til f_{j_1}}_{l_1}u_{\ify},
\label{tiltil}
\eeq
with the conditions
that $j_k\ne j_{k+1}$,
$\til e_{j_{k+1}} (\til f_{j_k}^m\cd u_{\ify})=0$
if $l'_{m-1}<k\leq l'_m$ and
$$
j_{l(b)}=
\left\{
\begin{array}{ll}
1&{\hbox{ if all walls in $b$ are }}+,\\
0&{\hbox{ if all walls in $b$ are }}-.
\end{array}
\right.
$$
\end{pr5}

To show this proposition, we shall see the following lemmas:

\begin{lem5}
\label{aik}
For $b=(\cd, i_k,i_{k+1},\cd, i_{-1})\in {\cal P}(\ify)$
and $i=0,1$
we set
\beq
A^{(i)}_k(b):=\sigma(i)(i_k+2\sum_{j<k}i_j),\qq
{\hbox{ where }}\sigma(i)=
\left\{
\begin{array}{ll}
+&{\hbox{ if }}i=1,\\
-&{\hbox{ if }}i=0.
\end{array}
\right.
\label{ak}
\eeq
\begin{enumerate}
\item
If there exists $k$ such that
$$
A^{(i)}_k(b)>A^{(i)}_{\nu}(b){\hbox{ for }}
  \nu<k{\hbox{ and }}
A^{(i)}_k(b)\leq A^{(i)}_{\nu}(b){\hbox{ for }} k<\nu\leq -1,
$$
we have
$\eit b=(\cd ,i_{k-1},\eit(i_k),i_{k+1},\cd,i_{-1}),$
otherwise, $\eit b=0$.
\item
If there exists $k$ such that
$$
A^{(i)}_k(b)\geq A^{(i)}_{\nu}(b){\hbox{ for }}
 \nu<k{\hbox{ and }}
A^{(i)}_k(b)> A^{(i)}_{\nu}(b){\hbox{ for }} k<\nu\leq -1,
$$
we have
$\fit b=(\cd ,i_{k-1},\fit(i_k),i_{k+1},\cd,i_{-1}).$
\item
$\vep_i(b)=\mathop{{\rm max}}_{k\leq -1}\{A^{(i)}_k(b)\}.$
\item
If $i_k=-i_{k+1}$, then $A^{(i)}_k(b)=A^{(i)}_{k+1}(b)$.
\item
If $\eit b=(\cd,i_{k-1},\eit(i_k),i_{k+1},\cd)$, we have
$$
A^{(i)}_{\nu}(\eit b)=
\left\{
\begin{array}{ll}
A^{(i)}_{\nu}(b)&{\rm if }\q \nu < k,\\
A^{(i)}_{\nu}(b)-1&{\rm if }\q \nu = k,\\
A^{(i)}_{\nu}(b)-2&{\rm if }\q \nu > k.
\end{array}
\right.
$$
\end{enumerate}
\end{lem5}

\nd
{\sl Remark.\,\,}
Since $A^{(i)}_j(b)=A^{(i)}_{j+1}(b)$ for $|j|\gg 0$,
there always exists $k$ as in (ii).
Thus,  $\fit b\ne 0$ for any $i$.

\vskip5pt
\nd
{\sl Proof.\,}
(i) and (ii) are easily obtained by Lemma 3.6 in \cite{K3}.
Using (\ref{tensor-vep}) repeatedly, we get (iii).
(iv) and (v) are trivial by the definition of $A^{(i)}_k(b)$.\qed

\begin{lem5}
\label{shift}
For $b=(\cd, i_{k},i_{k+1},\cd,i_{-1})
\in {\cal P}(\ify)\cong B(\ify)$
suppose that
all walls in $b$ are $-$ $($resp. $+$$)$
and set ${\til e_0}^{\vep_0(b)} b=(\cd,i'_{k},i'_{k+1},\cd,i'_{-1})$.
$($resp. ${\til e_1}^{\vep_1(b)}
b=(\cd,i'_{k},i'_{k+1},\cd,i'_{-1})$$)$.
 Then we have
\beq
i'_k=-i_{k-1}\qq{\hbox{ for }}k=-1,-2,\cd
\label{i=-i}
\eeq
\end{lem5}
{\sl Proof.\,}
We shall show the case that all walls are $-$.
Let  $M:=\{m_1,\cd, m_p\}$ be the set of
indices such that $i_{m_j-1}+i_{m_j}<0$ and $m_1<m_2<\cd <m_p$.
Here note that
any domain $i_{m_j},i_{m_j+1},\cd,i_{m_{j+1}-1}$
is in the form;
$l,-l,\cd,(-)^{m_{j+1}-m_j-1}l$.
By this fact and Lemma \ref{aik} (iii), we have
\beqn
&&A^{(0)}_{m_1}(b)<A^{(0)}_{m_2}(b)<\cd<A^{(0)}_{m_p}(b)=\vep_0(b),
\label{AAAA}\\
&& A^{(0)}_{m_j}(b)-A^{(0)}_{m_{j-1}}(b)=-i_{m_j-1}-i_{m_j}>0.
\eeqn
We set $n_j:=-i_{m_j-1}-i_{m_j}$.
Here note that by Lemma \ref{aik} (i) and (iv),
$\eit$ acts on some $i_k$ ($k\in M$) or $\eit b=0$.
In this case, by (\ref{AAAA}) we know that
$\til e_0$ acts on $i_{m_p}$.
Since $\til e_0(i_k)=i_k+1$, if we apply ${\til e_0}^{n_p}$
on $b$, by Lemma \ref{aik} (v) we get
\beqn
{\til e_0}^{n_p}b
&=& (\cd,i_{m_p-1},{\til e_0}^{n_p}(i_{m_p}),i_{m_p+1},\cd)\\
&=& (\cd,i_{m_p-1},-i_{m_p-1},i_{m_p+1},\cd)
\eeqn
Furthermore, if we apply ${\til e_0}$ on ${\til e_0}^{n_p}b$,
it acts on $i_{m_{p-1}}$ by Lemma \ref{aik} (ii) (v).
Repeating this,
we obtain
$$
{\til e_0}^{\vep_0(b)}b
=(\cd, {\til e_0}^{n_1}(i_{m_1}),\cd,{\til e_0}^{n_j}(i_{m_j}),
\cd,{\til e_0}^{n_p}(i_{m_p}),\cd).
$$
Here note that
$$
\sum_{j=1}^pn_j=-\sum_{j=1}^p(i_{m_j-1}+i_{m_j})
=A^{(0)}_{m_p}(b)=\vep_0(b).
$$
Since $i_{m_j-1}+{\til e_0}^{n_j}(i_{m_j})=0$, we have
\beq
i'_{m_j}=-i_{m_j-1}.
\label{im=-im}
\eeq
If $k\not\in M$, then $i'_k=i_k$ and $i_{k-1}+i_k=0$.
This implies
\beq
i'_k=-i_{k-1}\qq{\rm for }\q k\not\in M.
\label{ik=-ik}
\eeq
By (\ref{im=-im}) and (\ref{ik=-ik}) we obtained the desired result.
The case that all walls are $+$ is also shown similarly.
\qed

\vskip5pt\nd
{\sl Proof of Proposition \ref{b1}.\,}
We shall show by the induction on $l(b)$.
If $l(b)=1$, $b$ is in the following form;
$$
b=(\cd ,0,0,0,k).
$$
If $k> 0$, all walls in $b$ are $+$ and $b=\til f_1^ku_{\ify}$.
If $k< 0$, all walls in $b$ are $-$ and $b=\til f_0^{|k|}u_{\ify}$.
Suppose that $l(b)>1$ and all walls in $b$ are $-$.
%Let $d_1,d_2,\cd,d_p$ and
%$d'_1,d'_2,\cd,d'_{r}$ be the
%the finite domains with non-zero length in $b$ and
%${\til e_0}^{\vep_0(b)}b$ respectively.
%Here note that
%\beq
%{\hbox{if $l(d_p)>1$, $p=r$ and if $l(d_p)=1$, $r=p-1$.}}
%\label{prpr}
%\eeq
%Then by Lemma \ref{shift} and (\ref{prpr}), we have
%\beqn
%&&
%l(d'_j)=l(d_j)\q{\rm for }\q j=1,\cd, p-1,\\
%&&
%{\hbox{ if }}l(d_p)>1, l(d'_p)=l(d_p)-1.
%\eeqn
%where $l(d)$ is the length of $d$.
%%%%%%%%%
\nd
Lemma \ref{shift} implies that
all walls in ${\til e_0}^{\vep_0(b)}b$ are $+$ and
$l({\til e_0}^{\vep_0(b)}b)=l(b)-1$.
Then, by the hypothesis of the  induction,
\beq
{\til e_0}^{\vep_0(b)}b=
\underbrace{{\til f_1}^{n}
\cd{\til f_{j_{l'_{n-1}+1}}}^{n}}_{-1+l_{n}}
\cd\cd
\underbrace{\til f_{j_{l'_2}}^{2}
\cd\til f_{j_{l'_1+1}}^{2}}_{l_{2}}
\underbrace{\til f_{j_{l'_1}}
\cd\til f_{j_1}}_{l_1}u_{\ify}.
\label{tiltil}
\eeq
Since $\vep_0(b)=n$ and $b={\til f_0}^n{\til e_0}^nb$,
we obtain the desired result.
The case that all walls in $b$ are $+$
is also shown similarly.
\qed

\vskip5pt
\nd
We shall introduce the following lemma similar
 to Lemma \ref{aik}.

\begin{lem5}
\label{hat-a}
For $b=(\cd,a_k,a_{k-1},\cd,a_2,a_1)\in B_I$
we set
\beqn
\what A^{(0)}_{k}(b)
&:=&a_{2k-1}+2\sum_{j>k}(a_{2j-1}-a_{2j-2})\q(k\geq1),
\label{hat-a-0}\\
\what A^{(1)}_{k}(b)
&:=&a_{2k}+2\sum_{j>k}(a_{2j}-a_{2j-1})\q(k\geq1).
\label{hat-a-1}
\eeqn
\begin{enumerate}
\item
If there exists $k$ such that
$\what A^{(i)}_k(b)>\what A^{(i)}_{\nu}(b)$ for $k<\nu$ and
$\what A^{(i)}_k(b)\geq \what A^{(i)}_{\nu}(b)$ for $k>\nu\geq1$,
$$
\til e_ib=
\left\{
\begin{array}{ll}
(\cd, a_{2k-1}-1,\cd)&{\hbox{ if }}i=0,\\
(\cd,a_{2k}-1,\cd)&{\hbox{ if }}i=1,
\end{array}
\right.
$$
otherwise, $\eit b=0$.

\nd
If there exists $k$ such that
$\what A^{(i)}_k(b)\geq\what A^{(i)}_{\nu}(b)$ for $k<\nu$ and
$\what A^{(i)}_k(b)> \what A^{(i)}_{\nu}(b)$ for $k>\nu\geq 1$,
$$
\til f_ib=
\left\{
\begin{array}{ll}
(\cd, a_{2k-1}+1,\cd)&{\hbox{ if }}i=0,\\
(\cd,a_{2k}+1,\cd)&{\hbox{ if }}i=1,
\end{array}
\right.
$$
\item
$\vep_i(b)=\mathop{{\rm max}}_{k\geq 1}\{\what A^{(i)}_k\}.$
\end{enumerate}
\end{lem5}
It is easy to show this lemma by Lemma 1.3.6 in \cite{K3}.

%Let $b_1$ be an element as in Proposition \ref{b1}.
%Since $\vep_i^*({b_1}^*)=\vep_i({b_1})$ and by the fact that
%$\(b)=b'\ot {\fit}^{\vep_i^*(b)}(0)_i$,
%we have that
%\beqn
%\Psi(b_1)=u_{\ify}\ot&\cd &\ot
%\underbrace{\til f_{j_1}\ot\cd\ot
%\til f_{j_{l_1}}}_{l_1}
%\underbrace{{\til f_{j_{l_1+1}}}^2\ot\cd\ot
%{\til f_{j_{l'_2}}}^2}_{l'_2}\label{51}\\
%&&\ot\cd\cd\ot
%\underbrace{{\til f_{j_{l'_{n-1}+1}}}^n\ot\cd\ot
%{\til f_{j_{l_{l'_n}}}}^n}_{l'_n}u_{\ify},\nonumber
%\eeqn
%where we omit $(0)_j$.

\vskip10pt
We set
\beqn
\hspace{-10pt}B_{I}& := &\left\{(\cdot\cdot,a_k,a_{k-1},\cd,a_1)
|\hspace{-12pt}
\begin{array}{ll}
& a_k\in\ZZ_{\geq 0},\,a_{k+1} \leq a_{k}\\
&a_k=0\q{\rm for}\q k\gg 0
\end{array}
\right\}\subset{\Psi(B(\ify))} ,\\
\hspace{-10pt}B_{II}& := &
\left\{(\cd,i_k,i_{k+1},\cd,i_{-1})\in{\cal P}(\ify)\,|\hspace{-12pt}
\begin{array}{ll}
&i_{2k}\geq0, i_{2k-1}\leq 0,\\
&|i_{k-1}|\leq |i_{k}|.
\end{array}
\right\}.
\eeqn
where we set $(\cd,a_k,\cd,a_2,a_1)
:=u_{\ify}\ot\cd\ot {\til f_1}^{a_2}(0)_1\ot{\til f_0}^{a_1}(0)_0$.

\nd
The following lemma guarantees that $B_I$ and $B_{II}$ are
stable by the action of $\eit$.
\begin{lem5}
\label{stab12}
\begin{enumerate}
\item
If $b=(\cd,a_{k},a_{k-1},\cd,a_1)\in B_I$
and $\eit b\ne0$,
$\eit b\in B_I$.
\item
If $b=(\cd,i_k,i_{k+1},\cd,i_{-1})\in B_{II}$ and $\eit b\ne0$,
$\eit b\in B_{II}$.
\end{enumerate}
\end{lem5}

\nd
{\sl Proof.\,\,}(i)\,\,\,
We consider the case of $i=0$ and $b\in B_I$.
We assume that
$a_{2k+1}$ changes to $a_{2k+1}-1$  by the action of $\til e_0$
and $a_{2k+1}=a_{2k+2}$.
$\what A^{(0)}_{k}(b)-\what A^{(0)}_{k+1}(b)=
a_{2k+1}-2a_{2k+2}+a_{2k+3}=a_{2k+3}-a_{2k+2}\leq 0$.
This means $\what A^{(0)}_{k}(b)\leq \what A^{(0)}_{k+1}(b)$,
which contradicts Lemma \ref{hat-a}.
Thus
if  $a_{2k+1}$ changes to $a_{2k+1}-1$ by the action of $\til e_0$,
$a_{2k+2}<a_{2k+1}$.
As for the $i=1$-case, we can show similarly.

\nd
(ii)\,\,\,
For $b=(\cd,i_k,i_{k+1},\cd,i_{-1})\in B_{II}$
suppose that $\til e_0b=(\cd,i_k+1,\cd)\ne0$.
Let $A^{(i)}_k(b)$ be as in Lemma \ref{aik}.
By the assumption we have $A^{(0)}_k(b)>A^{(0)}_{k-1}(b)$ and then
$i_k+i_{k-1}<0$. Suppose that $i_k\geq 0$. Then
 $-i_{k-1}>i_k\geq 0$.
This means $|i_k|<|i_{k-1}|$, which is a contradiction.
Thus we have $i_k<0$. This implies $0\leq i_{k-1}< -i_k$
and then $0\leq i_{k-1}\leq -i_k-1$. Therefore, we get
$|i_{k-1}|\leq |i_k+1|=|\til e_0(i_k)|$.
Then it follows that if $\til e_0b\ne 0$, $\til e_0b\in B_{II}$.
The $i=1$-case can be shown by the similar way.\qed

\vskip5pt
\nd
{\sl Remark.\,}
\begin{enumerate}
\item
We know that a vector in $B_I$ is in the image of
$\Psi$ by Proposition \ref{b-inf-sl2}.
\item
Any element $b=(\cd,a_2,a_1)\in B_I$ is in the form:
$$
b=(\cd,0,0,\underbrace{1\cd,1}_{l_1},\underbrace{2,\cd,2}_{l_2},
\cd,\underbrace{n,\cd,n}_{l_n}).
$$
%Then for an element in $B_I$ we call a maximal subsequence
%$d_m=m,m,\cd,m$ a {\it domain} and
%$l_m$ the {\it length} of $d_m$.
\end{enumerate}

Now, we shall show the following proposition:
\begin{pr5}
\label{taiou}
For $b\in B(\ify)$, suppose that
\beq
\Psi(b) =(  \cd 0,0, d_1,\cd, d_n)\in B_I,\label{52}
\eeq
where as in Remark (ii) as above;
$$
d_m = \underbrace{m,\cd,m}_{l_m}=
\underbrace{{\til f_{s^m_1}}^m\ot\cd\ot {\til f_{s^m_{l_m}}}^m}_{l_m}
\in B_{s^m_1}\ot\cd\ot B_{s^m_{l_m}},\q(s^m_j=0,1).
$$
Then we have
\beq
\Phi(b) =(  \cd 0,0,\hat d_1,\cd,\hat d_n)\in B_{II},\label{53}
\eeq
where $\Phi$ is the isomorphism
$\Phi:B(\ify)\mapright{\sim}{\cal P}(\ify)$ and
$$
\hat d_m=
\underbrace{\sigma(s^m_1)m,\cd,\sigma(s^m_{l_m})m}_{l_m}\in
{B_{\ify}}^{\ot l_m},
$$
\end{pr5}

\vskip5pt
\nd
To show the proposition, we shall prepare several lemmas.

\begin{lem5}
\label{A=A}
Let $A^{(i)}_k(b)$ be as in Lemma \ref{aik} and $\what A^{(i)}_k(b)$
be as in Lemma \ref{hat-a}.
Let $b_1=(\cd,a_2,a_1)\in B_I$ be as in (\ref{52}) and
$b_2=(\cd, i_{-2},i_{-1})\in B_{II}$ be as in (\ref{53}).
Then
$\fit b_1=(\cd,a_k+1,\cd)$ $($resp. $\eit b_1=(\cd,a_k-1,\cd)$$)$
if and only if
$\fit b_2=(\cd,\fit(i_{-k}),\cd)$
$($resp. $\eit b_2=(\cd,\eit(i_{-k}),\cd)$$)$.
\end{lem5}
{\sl Proof.\,}
First, we shall show
\beq
\what A^{(i)}_{k}(b_1)= A^{(i)}_{-2k+1-i}(b_2),
\label{a=hat-a}
\eeq
for any $k\in \ZZ_{>0}$ and $i\in\{0,1\}$.
For this purpose, we shall see
the following;

\begin{lem5}
\label{numbers}
\begin{description}
\item[(i)]
%For $b$ as in (\ref{52}) and $i\in\{0,1\}$,
If $a_{2k-1+i}$ is in $d_m$
({\it i.e.} $a_{2k-1+i}=m$),
\beq
\what A^{(i)}_{k}(b_1)=
\sharp\{j|s^j_1=i,\,\,1\leq j\leq m\}
-\sharp\{j|s^j_1=1-i,\,\,1\leq j\leq m\}.
\label{a=1-0}
\eeq
\item[(ii)]
%For $b$ as in (\ref{53}) and $i\in\{0,1\}$,
If $i_{-2k+1-i}$ is in $\hat d_m$,
\beq
 A^{(i)}_{-2k+1-i}(b_2)
=\sharp\{j|s^j_1=i,\,\,1\leq j\leq m\}
-\sharp\{j|s^j_1=1-i,\,\,1\leq j\leq m\}.
\label{hat-a=1-0}
\eeq
\end{description}
\end{lem5}
\vskip5pt
\nd
{\sl Remark.\,\,}
Even if $l_m=0$, $s^m_1$ and $s^m_{l_m}$ are uniquely determined
 by applying
the conditions $s^m_{l_m}\ne s^{m+1}_1$,
$s^m_1\ne s^{m-1}_{l_{m-1}}$ and if $l_m$ is even, $s^m_1\ne s^m_{l_m}$
to the nearest non-empty domains.  That is,
let  $d_a,d_b,d_c$ ($a<b<c$) be domains such that
$d_a$ and $d_c$ are non-empty and
$d_b$ is the empty domain surrounded by $d_a$ and $d_c$.
Then $s^b_1$ and $s^b_{l_b}$ are given by
$s^b_1=1-s^a_{l_a}$ and $s^b_{l_b}=1-s^c_1$.

\vskip6pt
\nd
{\sl Proof.\,\,}
(i)\,\,\,
We shall consider the $i=0$ case.
We know that
$a_{2k-2}$ is the left-most entry of some $d_m$ if and
only if $a_{2k-1}<a_{2k-2}$.
Now, $a_{2k-2}$ implies ${\til f_1}^{a_{2k-2}}(0)_1$.
Since the index of the left-most entry
of $d_j$ is $s^j_1$,
if $a_{2k-1}$ is in some $d_m$ ({\it i.e.} $a_{2k-1}=m$),
by the remark as above
 and (\ref{hat-a-0}), we have
\beqnn
\what A^{(0)}_k(b_1)
             &=& a_{2k-1}-2\sharp\{j|s^j_1=1,\;\;1\leq j\leq m\}\\
             &=& m- 2\sharp\{j|s^j_1=1,\,\,1\leq j\leq m\}\\
             &=& \sharp\{j|s^j_1=0,\,\,1\leq j\leq m\}
+\sharp\{j|s^j_1=1,\,\,1\leq j\leq m\}
-2\sharp\{j|s^j_1=1,\,\,1\leq j\leq m\}\\
         &=&\sharp\{j|s^j_1=0,\,\,1\leq j\leq m\}
           -\sharp\{j|s^j_1=1,\,\,1\leq j\leq m\}
\eeqnn
$\what A^{(1)}_k(b_1)$ is also given similarly.

\nd
(ii)\,\,\,
We shall consider the $i=0$ case.
By (\ref{ak}), we can write
\beq
A^{(0)}_k(b_2)=-\sum_{j\leq k}(i_{j-1}+i_j).
\label{aii}
\eeq
We know that $i_{j-1}+i_j\ne0$ if and only if
$i_j$ is the left-most element of some $\hat d_m$.
Here note that if $j$ is odd, $i_{j-1}+i_j\leq 0$ and
if $j$ is even, $i_{j-1}+i_j\geq 0$ by the conditions of $B_{II}$.
Let $i_j$ be the left-most element of $\hat d_m$.
If $j$ is even, $s^m_1=1$ and if $j$ is odd, $s^m_1=0$.
Thus, if $i_{-2k+1}$ is in $\hat d_m$, by
the remark as above and
(\ref{aii}) we have
\beqn
A^{(0)}_{-2k+1}(b_2) &=&-(\sharp\{j|s^j_1=1,\,\,1\leq j\leq m\}
   -\sharp\{j|s^j_1=0,\,\,1\leq j\leq m\})\\
  &=& \sharp\{j|s^j_1=0,\,\,1\leq j\leq m\}
     -\sharp\{j|s^j_1=1,\,\,1\leq j\leq m\}.
\eeqn
$A^{(1)}_{-2k}(b_2)$ is also obtained by the simlar argument.
\qed

\vskip5pt
\nd
{\sl Proof of lemma \ref{A=A}.\,\,}
By Lemma \ref{numbers} we obtain (\ref{a=hat-a}).

For $i=0$ case, we set
$\til f_0(b_1)=(\cd,a_{2k-1}+1,\cd)$.
By Lemma \ref{hat-a}
this implies that
$\what A^{(0)}_{k}(b_1)\geq \what A^{(0)}_{\nu}(b_1)$
for $k<\nu$ and
$\what A^{(0)}_{k}(b_1)> \what A^{(0)}_{\nu}(b_1)$
for $k>\nu\geq 1$.
Thus, by (\ref{a=hat-a}) we have
$A^{(0)}_{-2k+1}(b_2)\geq A^{(0)}_{-2\nu+1}(b_2)$
for $k<\nu$ and
$A^{(0)}_{-2k+1}(b_2)> A^{(0)}_{-2\nu+1}(b_2)$
for $k>\nu\geq1$.
Furthermore, by simple calculations and
the definition of $B_{II}$, for any $j\in \ZZ_{>0}$ we get
\beqn
&&
A^{(0)}_{-2j}(b_2)-A^{(0)}_{-2j+1}(b_2)
=i_{-2j+1}+i_{-2j}\leq 0,\\
&&
A^{(0)}_{-2j-1}(b_2)-A^{(0)}_{-2j}(b_2)
=i_{-2j-1}+i_{-2j}\geq 0.
\eeqn
This implies that
$$
A^{(0)}_{-2j-1}(b_2)\geq
A^{(0)}_{-2j}(b_2)\leq
A^{(0)}_{-2j+1}(b_2).
$$
Then by Lemma \ref{aik} (ii), we know that
the action by $\til f_0$ on $b_2$ never touch $i_{-2j}$.
Then we get
$\til f_0 b_2=(\cd,\til f_0(i_{-2k+1}),\cd)$.
Arguing similarly, we have that
$\til f_0 b_1=(\cd, a_{2k+1}+1,\cd)$ if
$\til f_0 b_2=(\cd,\til f_0(i_{-2k+1}),\cd)$.
The other cases are shown  similarly.
\qed

\begin{lem5}
\label{l=l-1}
For $b\in B_I$ $($resp. $b\in B_{II}$$)$
there exists $i\in\{0,1\}$ such that
\beq
l(\eit^{\vep_i(b)}(b))=l(b)-1,
\label{e=e-1}
\eeq
where $l(b)$ is a length of $b$ given as
the largest number $k>0$ such that
$a_k\ne0$
for $b=(\cd,a_j,a_{j-1},\cd,a_1)$
$($resp. $i_{-k}\ne(0)$ for $b=(\cd,i_{-k},i_{-k+1},\cd,i_{-1})$$)$.
\end{lem5}

\nd
{\sl Proof of Lemma \ref{l=l-1}.\,}
If $l(b)$ is even, we choose $i=1$ and if
$l(b)$ is odd, we choose $i=0$.
(Here note that we chose the sequence of indices :$(\cd,1,0,1,0)$
for $B_I$.).
We assume $l(b)$ is odd and then $i=0$, and set $l(b)=2k-1$ and
$$
{\til e_0}^{\vep_0(b)}b=(\cd,b_{2k},b_{2k-1},\cd,b_1)
\q{\hbox{(resp. }}(\cd,j_{-2k},j_{-2k+1},\cd,j_{-1})).
$$
Suppose that $b_{2k-1}>0$ (resp. $j_{-2k+1}>0$). Then we have
$\what A^{(0)}_k({\til e_0}^{\vep_0(b)}b)=b_{2k-1}>0$
(resp. $A^{(0)}_{-2k+1}({\til e_0}^{\vep_0(b)}b)=j_{-2k+1}>0$)
and then by Lemma \ref{hat-a} (ii) (resp. Lemma \ref{aik}(iii))
$\vep_0({\til e_0}^{\vep_0(b)}b)>0$, which is a contradiction
(see Example \ref{Example:crystal} (iv)).
Thus, we get $b_{2k-1}=0$ (resp. $j_{-2k+1}=0$).
For  $b\in B_I$, $a_{2j}$ implies ${\til f_1}^{a_{2j}}(0)_1$.
Then $\til e_0$ never touch $a_{2k-2}$. Then we have $a_{2k-2}=b_{2k-2}$.
By the definition of $B_{I}$, we know that
$b_{2k-2}=a_{2k-2}\geq a_{2k-1}>0$. Now, we obtain (\ref{e=e-1})
for $b\in B_I$.

\nd
For $b\in B_{II}$, by the assumption
$l(b)$ is odd and by the definition of $B_{II}$, we have
$i_{-2k+1}+i_{-2k+2}\geq 0$.
Then, by the facts that $\vp_0(i_{-2k+1})=i_{-2k+1}$
and $\vep_0(i_{-2k+2})=-i_{-2k+2}$,
we get $\vp_0(i_{-2k+1})\geq \vep_0(i_{-2k+2})$
and in general,
$\vp_0({\til e_0}^m(i_{-2k+1}))=\vp_0(i_{-2k+1})+m
\geq \vep_0(i_{-2k+2})=-i_{-2k+2}$ for $m\geq 0$.
Thus, by (\ref{tensor-e}) we obtain
$$
{\til e_0}^m(i_{-2k+1}\ot i_{-2k+2})
={\til e_0}^m(i_{-2k+1})\ot i_{-2k+2}.
$$
This implies that
action of ${\til e_0}^m$ ($0\leq m\leq \vep_0(b)$) never touch
the entry $i_{-2k+2}$ and then $i_{-2k+2}=j_{-2k+2}\ne0$.
Now, we get (\ref{e=e-1}) for $b\in B_{II}$.
\qed

\vskip5pt
\nd
{\sl Proof of Proposition \ref{taiou}.\,}
Let us complete the proof of Proposition \ref{taiou}
by using the induction on $l(\Psi(b))$.
For the case $l(\Psi(b))=1$,
we can write
$\Psi(b)=(\cd,0,0,a)=u_{\ify}\ot {\til f_0}^a(0)_0$ ($a>0$).
Then we obtain $\Phi(b)={\til f_0}^a(\cd,0,0,0)=
(\cd,0,0,-a)\in B_{II}$.
We assume $l(\Psi(b))>1$.
Set $\Psi(b)=(\cd,a_k,a_{k-1},\cd,a_1)$ and
$\Phi(b)=(\cd,i_k,i_{k+1},\cd,i_{-1})$.
By Lemma \ref{l=l-1}, there exists $i$ such that
$l(\eit^{\vep_i(b)}\Psi(b))=l(\Psi(b))-1$.
Since $B_I$ is stable by the action of $\eit$ by
Lemma \ref{stab12},  we can set
$\eit^{\vep_i(b)}\Psi(b)=(\cd,a'_k,\cd,a'_2,a'_1)=
(\cd, 0,0,d'_1, d'_2,\cd, d'_r)\in B_I$
with
$$
 d'_m=\underbrace{m,\cd,m}_{l'_m}={\til f_{t^m_1}}^m\ot\cd\ot
{\til f_{t^m_{l'_m}}}^m,
$$
where $t^m_j\in \{0,1\}$ and $l'_m$ is the length of $d'_m$.
By the hypothesis of the induction, we can write
$\eit^{\vep_i(b)}(\Phi(b))=\Phi(\eit^{\vep_i(b)} b)=
(\cd,i'_{-k},\cd,i'_{-2},i'_{-1})=
(\cd 0,0,\hat d'_1,\hat d'_2,\cd,\hat d'_r)$
with
$$
\hat d'_m=\sigma(t^m_1)m,\cd,\sigma(t^m_{l'_m})m.
$$
We consider $i=0$ case.
If $\til f_0{\til e_0}^{\vep_0(b)}\Psi(b)
=(\cd,a'_{2k+1}+1,\cd)$, by Lemma \ref{A=A}
we get
$\til f_0{\til e_0}^{\vep_0(b)}\Phi(b)
=(\cd,\til f_0(i'_{-2k-1}),\cd)$.
Here note that
$\til f_0(i'_{-2k-1})=i'_{-2k-1}-1$ and $i'_{-2k-1}\leq 0$
by the definition of $B_{II}$, then
$|\til f_0(i'_{-2k-1})|=|i'_{-2k-1}|+1=a'_{2k+1}+1$.
Repeating this $\vep_0(b)$ times,
since ${\til f_0}^{\vep_0(b)}{\til e_0}^{\vep_0(b)}\Phi(b)=\Phi(b)$
and
${\til f_0}^m{\til e_0}^{\vep_0(b)}\Psi(b)\in B_I$
($0\leq m\leq \vep_0(b)$),
we obtain the desired result.
\qed

\begin{pr5}
\label{psi-phi}
For $b\in B(\ify)$ suppose that
all walls in
$\Phi(b)=$\hfill\break
$(\cd,i_k,i_{k+1},\cd,i_{-1})\in {\cal P}(\ify)$
are $-$ $($resp. $+$$)$ and the total number of walls is $n\in\ZZ_{>0}$.
$($Then there are $n$ finite domains.$)$.
Let $l_j$ $(1\leq j\leq n)$ be the length of the $j$-th finite
domain. Then we have
$$
\Phi(b^*)=(\cd,0,0,\bar d_1,\bar d_2,\cd,\bar d_n),
$$
with
\beqn
&&\bar d_j
=\underbrace{
(-)^{k_j}j,(-)^{k_j+1}j,\cd,
(-)^{k_j+l_j-1}j}_{l_j},
\label{desc-star}\\
&&({\rm resp. }\,\,
\bar d_j
=\underbrace{
(-)^{k_j+1}j,(-)^{k_j+2}j,\cd,
(-)^{k_j+l_j}j}_{l_j}),
\eeqn
where $k_j\in \{-1,-2,\cd\}$ is the position of the left-most
entry of $\bar d_j$.
\end{pr5}

\vskip5pt
\nd
{\sl Proof.\,}
We assume that all walls in $b$ are $-$.
Here note the following:

\nd
if $b\in B(\ify)$ can be written $b=\fit^mb'$ with
$\eit b'=0$ (that is, $\vep_i(b)=m$),
\beq
\Psi_i(b^*)=b'\ot \fit^m(0)_i
\label{psi-star}
\eeq
since $\vep^*_i(b^*)=\vep_i(b)=m$. (See 5.1.).
By Proposition \ref{b1}, we can write $b$ in the following form;
$$
b=
\underbrace{\til f_0^{n}
\cd\til f_{j_{l'_{n-1}+1}}^{n}}_{l_{n}}
\cd\cd
\underbrace{\til f_{j_{l'_2}}^{2}
\cd\til f_{j_{l'_1+1}}^{2}}_{l_{2}}
\underbrace{\til f_{j_{l'_1}}
\cd\til f_{j_1}}_{l_1}u_{\ify},\q(l'_m=\sum_{j=1}^ml_j)
$$
with the conditions
$j_k\ne j_{k+1}$ ($1\leq k<l'_n$) and
\beq
\til e_{j_{k+1}} (\til f_{j_k}^m\cd u_{\ify})=0,
\label{eee}
\eeq
for $l'_{m-1}<k\leq l'_m$ and $1\leq m\leq n$.
By virtue of (\ref{eee}), we can apply (\ref{psi-star})
to $b^*$  repeatedly and  obtain;
\beq
\Psi(b^*)=
(\cd ,0,0,d_1,d_2,\cd,d_n),
\label{b-star}
\eeq
where
$$
d_m=\underbrace{m,\cd,m}_{l_m}=
\underbrace{
{\til f_{s^m_1}}^m\ot\cd\ot{\til f_{s^m_{l_m}}}^m}_{l_m}
\q{\rm and}\q s^n_{l_n}=0\q(1\leq m\leq n).
$$
Since the vector (\ref{b-star}) is in $B_I$,
 by Proposition \ref{taiou}
we have
$$
\Phi(b^*)=(\cd,0,0,\hat d_1,\hat d_2,\cd,\hat d_n),
$$
with
$$
\hat d_m=\underbrace{
\sigma(s^m_1)m,\cd,\sigma(s^m_{l_m})m}_{l_m}\in B_{\ify}^{\ot l_m}.
$$
We set
$\Phi(b^*)=(\cd,i_{-k},\cd,i_{-1})$
and  $\Psi(b^*)=(\cd,a_{k},\cd,a_1)$.
If $k$ is odd (resp. even),
we have $a_k={\til f_0}^{a_k}(0)_0$
(resp. $a_k={\til f_1}^{a_k}(0)_1$)
and $i_{-k}<0$
(resp. $i_{-k}>0$).
Then we have $\hat d_j=\bar d_j$ ($1\leq j\leq n$).
\qed

\vskip10pt
Now, we shall give the similar description for $B(-\ify)$.
By the definition of $\wedge$, we have
\beqnn
B_{\ify}^{\wedge}&\mapright{\sim}&B_{\ify}\\
     (n)^{\wedge}&\mapsto &(-n).
\eeqnn
As mentioned in Example \ref{Example:crystal} (iv),
$B(\mp\ify)=B(\pm\ify)^{\vee}\cong B(\pm\ify)^{\wedge}$.
Then we can identify $\vee$ with $\wedge$ on $B(\pm\ify)$.
We obtain the following isomorphism:
$$
\Phi^+:=\wedge\circ\Phi\circ\wedge:
B(-\ify)\mapright{\sim}{\cal P}(-\ify).
$$

\vskip6pt
\nd
{\sl Remark.\,}
For $b\in B(\ify)$, if we set
\beqnn
&&\Phi(b)=(\cd, i_{-k},i_{-k+1},\cd,i_{-1}),\\
&&\Phi^+(b^{\wedge})=(j_0,\cd,j_k,j_{k+1},\cd),
\eeqnn
we get $j_k=-i_{-k-1}$.

\vskip6pt
\nd
\begin{pr5}
\label{u-plus}
For $b\in B(-\ify)$ suppose that
all walls in
$\Phi^+(b)=(i_0,\cd,i_k,i_{k+1},\cd)\in {\cal P}(-\ify)$
are $-$ $($resp. $+$$)$ and the number of walls is $n\in\ZZ_{>0}$.
(Then there are $n$ finite domains.).
Let $l_j$ $(1\leq j\leq n)$ be the length of the finite
domain $\bar d_j$. Then we have
$$
\Phi^+(b^*)=(\bar d_n,\bar d_{n-1},\cd,\bar d_1,0,0,\cd)
$$
where
\beqn
&&\bar d_j
:=\underbrace{
(-)^{k_j+1}j,(-)^{k_j+2}j,\cd,
(-)^{k_j+l_j}j}_{l_j},
\label{desc-star}\\
({\rm resp.}&&
\bar d_j
:=\underbrace{
(-)^{k_j}j,(-)^{k_j+1}j,\cd,
(-)^{k_j+l_j-1}j}_{l_j}),
\eeqn
where $k_j\in\{0,1,2,\cd\}$ is the position of the left-most
entry of $\bar d_j$.
\end{pr5}

\vskip5pt
\nd
{\sl Proof.\,\,}
If all the walls in some element
$b_0\in B(\ify)$ are $+$ (resp. $-$),
by the Remark as above
all the walls in $b:=b_0^{\wedge}=b_0^{\vee}$ are $-$ (resp. $+$).
By the commutativity of $\vee$ and $*$ and Proposition \ref{psi-phi},
 we get the desired result.\qed

\vskip6pt
\nd
For $m\in \ZZ$, we set
$$
{\cal P}_m(-\ify)=
\left\{(i_0,i_1,\cd,i_k,\cd)|
\hspace{-3pt}
\begin{array}{l}
i_k\in B_{\ify}\q{\rm and \,\,\,if }\q |k|\gg 0, \\
i_{-k}=(0), \q i_{2k}=(m), \q i_{2k+1}=(-m)
\end{array}\right\}.
$$
For $\lm=m(\Lm_0-\Lm_1)+l\del$  by \cite[5.2]{N},
there exists the isomorphism
\beqn
T_{\lm}\ot {\cal P}(-\ify)& \cong &{\cal P}_m(-\ify)
\label{TPP}\\
t_{\lm}\ot (i_0,\cd,i_{2k},i_{2k+1},\cd)&\leftrightarrow&
(i_0+m,\cd,i_k+(-)^km,\cd)\nonumber
\eeqn
where weight of $b$ in ${\cal P}_m(-\ify)$ is given by
$wt(b)=\lm+wt(b')$ ($t_{\lm}\ot b'\leftrightarrow b$).

\nd
Now, we can describe the operation $*$ on extremal vectors.
For a level 0 weight $\lm=m(\Lm_0-\Lm_1)+l\del$,
let $b=(\cd,i_k,i_{k+1},\cd)\in B(\util_0)_{-\lm}$
be an extremal vector.
By (\ref{m01}) and Lemma \ref{mmm},
if $m>0$ $($resp. $m<0$$)$, then all walls in $b$ are $-$
$($resp. $+$$)$ and the number of walls is $|m|$.
There are $|m|-1$ domains in $b$ denoted $d_1,d_2,\cd,d_{|m|-1}$.
Let $l_j$ be the length of $d_j$ and $k_j$
be the position of the left-most entry in
a domain $d_j$ with $l_j>0$.

\begin{thm5}
\label{star-on-ext}
For $\lm=m(\Lm_0-\Lm_1)+l\del$ with $m>0$ $($resp. $m<0$$)$
let $b=(\cd,i_k,i_{k+1},\cd)\in B(\util_0)_{-\lm}$
be an extremal vector as above.
Set
\beqn
&&\bar d_j:=i'_{k_j},i'_{k_j+1},\cd,i'_{k_{j+1}-1}
=\underbrace{(-)^{k_j}j,(-)^{k_j+1}j,\cd,(-)^{k_{j+1}-1}j}_{l_j}
\label{d-bar}\\
&&\hspace{-40pt}({\rm resp.}\q
\bar d_j:=i'_{k_j},i'_{k_j+1},\cd,i'_{k_{j+1}-1}
=\underbrace{(-)^{k_j+1}j,(-)^{k_j+2}j,\cd,(-)^{k_{j+1}}j}_{l_j}).
\eeqn
Then we obtain
$$
b^*=(\cd,i'_{k_j},i'_{k_j+1},\cd)
=(\cd,0,0,\bar d_1,\cd,\bar d_{|m|-1},(-)^{t+1}m,(-)^{t+2}m,\cd),
$$
where $t$ is the position of the right-most entry in the subsequence
$\bar d_1,\cd,\bar d_{|m|-1}$, that is,
$t=k_j+l_j+l_{j+1}+\cd+l_{|m|-1}-1$ if $l_j\ne 0$.
\end{thm5}

\vskip6pt
\nd
{\sl Proof of Theorem \ref{star-on-ext}.\,}
We assume $m>0$ and then all walls in $b$ are $-$.
Set $b=b_1\ot t_{\mu}\ot b_2
\in B(\ify)\ot T_{\mu}\ot B(-\ify)$ ($\mu:=-\lm-wt(b_1)-wt(b_2)$),
and
\beqnn
&&b_1=(\cd, i_k,\cd,i_{-2},i_{-1}),\\
&&t_{\mu}\ot b_2=(i_0,i_1,\cd,i_k,\cd).
\eeqnn
Let $m_1$ and $m_2$ be the total number of walls
in $b_1$ and $t_{\mu}\ot b_2$ respectively.
Here note that $m_1+m_2\leq m$ since
the walls of $b$ at position 0 are not included
in $b_1$ nor $t_{\mu}\ot b_2$.
Let $d^1_j$ ($1\leq j\leq m_1$) and $d^2_j$ ($1\leq j\leq m_2$)
be the domain in $b_1$ and $t_{\mu}\ot b_2$ respectively and
set $l^1_j$ and $l^2_j$ the length of $d^1_j$ and $d^2_j$ respectively
and
$k^1_j$ and $k^2_j$ the position of the left-most entry in
$d^1_j$ and $d^2_j$ with non-zero length respectively.
That is,
\beqnn
&&b_1=(\cd,0,0,d^1_1,d^1_2,\cd,d^1_{m_1}),\\
&&t_{\mu}\ot b_2=(d^2_{m_2},d^2_{m_2-1},\cd,d^2_1,0,0,\cd).
\eeqnn
Here note that
$l^1_j=l_j$ for $(1\leq j<m_1)$ and
$l^2_j=l_{m+1-j}$ for $(1\leq j< m_2)$.
Also note that
$k^1_j=k_j$ if $l^1_j>0$ ($1\leq j\leq m_1$)
 and $k^2_j=k_{m+1-j}$ if $l^2_j>0$ $(1\leq j<m_2)$.
By Proposition \ref{taiou}, Proposition \ref{u-plus} and (\ref{TPP}),
we have
\beqn
&&b_1^*=
(\cd,0,0,\hat d^1_1,\cd,\hat d^1_{m_1}),\\
&&t_{\lm}\ot b_2^*=
(\hat d^2_{m_2},\hat d^2_{m_2-1},
\cd,\hat d^2_1,(-)^{t+1}m,(-)^{t+2}m,\cd),
\eeqn
where $t$ is the position of the right-most entry in
the subsequence
$\hat d^2_{m_2},\hat d^2_{m_2-1},
\cd,\hat d^2_1$ and
\beqn
\hat d^1_j&=&
\underbrace{(-)^{k^1_j}j,(-)^{k^1_j+1}j,
\cd,(-)^{k^1_{j+1}-1}j}_{l^1_j},\\
\hat d^2_j&=&
\underbrace{(-)^{k^2_j}(m-j),(-)^{k^2_j+1}(m-j),
\cd,(-)^{k^2_{j-1}-1}(m-j)}_{l^2_j}.
\eeqn
In particular, if we denote
$b^*=b_1^*\ot t_{\lm}\ot b_2^*=(\cd,j_k,j_{k+1},\cd)$,
we have
\beq
j_{-1}=-m_1\q{\rm and }\q j_0=m-m_2.
\label{j-1j0}
\eeq
This implies there are
$j_{-1}+j_0=m-m_1-m_2\geq 0$ walls at position 0.
We know that there are $m_1$ walls in $b_1^*$ and
$m_2$ walls in $t_{\lm}\ot b_2^*$.
Thus, the total number of walls in $b^*$
is $(m_1+m_2)+(m-m_1-m_2)=m$ and there are $m-1$ finite domains in $b^*$.
We denote them $\hat d_j$ ($1\leq j\leq m-1$).
There are the following two cases:
\begin{description}
\item[(I)]
There is no wall at position 0 in $b$ (that is, $m_1+m_2=m$).
\item[(II)]
There are walls at position 0 in $b$ (that is, $m_1+m_2<m$).
\end{description}

\vskip5pt
\nd
(I)\,\,\,\,
In this case, $m_1=m-m_2$. Then by (\ref{j-1j0}) we know that
$j_{-1}+j_0=0$ and then
there is no wall at position 0 in $b^*$.
Thus, we have
\beqnn
&&\hat d_j=\hat d^1_j=\bar d_j\q(1\leq j< m_1)\q{\rm and }\q
\hat d_j=\hat d^2_{m-j}=\bar d_j\q(m-m_2=m_1<j\leq m-1),\\
&&\hat d_{m_1}=\hat d^1_{m_1}\cup \hat d^2_{m_2}
=(-)^{k_{m_1}}m_1,(-)^{k_{m_1}+1}m_1,\cd,(-)^{k_{m_1+1}-1}m_1
=\bar d_{m_1}.
\eeqnn

\nd
(II)\,\,\,\,
In this case, $m-m_1-m_2>0$. This menas that
there exist $m-m_1-m_2$ walls at positin 0 in $b^*$ by (\ref{j-1j0}).
\beqnn
&&\hat d_j=\hat d^1_j=\bar d_j\q(1\leq j< m_1)\q{\rm and }\q
\hat d_j=\hat d^2_{m-j}=\bar d_j\q(m-m_2<j\leq m-1),\\
&& \hat d_j=\emptyset=\bar d_j\q(m_1\leq j \leq m-m_2).
\eeqnn
Now, we obtain the desired result.
The case $m<0$ is also shown by the simialr way.\qed

\begin{ex5}
For
$b=(\cd,0,0,\underbrace{
\stackrel{-5}{-1},\stackrel{-4}{1}}_{d_1},
\underbrace{\stackrel{-3}{-2},
\stackrel{-2}{2},\stackrel{-1}{-2}}_{d_2},
\underbrace{\stackrel{0}{1},\stackrel{1}{-1},\stackrel{2}{1}}_{d_3},
\stackrel{3}{-2},\stackrel{4}{2},-2,2,\cd)$, we have
$$
b^*=(\cd,0,0,\underbrace{
\stackrel{-5}{-1},\stackrel{-4}{1}}_{\bar d_1},
\underbrace{\stackrel{-3}{-2},
\stackrel{-2}{2},\stackrel{-1}{-2}}_{\bar d_2},
\underbrace{\stackrel{0}{3},\stackrel{1}{-3},\stackrel{2}{3}}_{\bar d_3},
\stackrel{3}{-4},\stackrel{4}{4},-4,4,\cd).
$$
\end{ex5}

\vskip7pt
\nd
\newtheorem{cor5}[thm5]{Corollary}
\begin{cor5}
\label{b-max-final}
For level 0 weight $\lm=m(\Lm_0-\Lm_1)+l\del$ ($l,m\in\ZZ$),
\beq
\bmax(\lm)=\bigoplus_{\vec c\in \ZZ^{|m|-1}_{\geq 0}}
{\cal P}_{m,l}(|m|;\vec m;\vec c),
\label{bmax-star}
\eeq
where
$$
\vec m=
\left\{
\begin{array}{ll}
(1,2,\cd,m-1)&{\rm if }\q m\geq 0,\\
(-1,-2,\cd,m+1)&{\rm if }\q m<0.
\end{array}
\right.
$$
\end{cor5}

\vskip6pt
\nd
{\sl Proof.\,}
Let $b$ and $\{\bar d_j\}$ be as in Theorem \ref{star-on-ext}.
By (\ref{d-bar})
$t(\bar d_j)=$ the type of domain $\bar d_j$
(see \cite[Sec 7.1]{N})
is given by
$$
t(\bar d_j)
=\left\{
\begin{array}{ll}
j&\q{\rm if }\q m>0,\\
-j&\q{\rm if}\q m<0.
\end{array}
\right.
$$

For $b^*$ as in Theorem \ref{star-on-ext}
setting $c_j:=[[l_j/2]]$ ($[[n]]$ is the Gauss's symbol),
we get
\beq
b^*\in {\cal P}_m(|m|;\vec m;\vec c),
\eeq
where $\vec c=(c_1,c_2,\cd, c_{|m|-1})$. This means
$$
\bmax(\lm)\subset \bigoplus_{\vec c\in \ZZ^{|m|-1}_{\geq 0}}
{\cal P}_{m,l}(|m|;\vec m;\vec c).
$$
Now, we assume that $m>0$ (resp. $m<0$).
Let $b_0$ be an extremal vector in ${\cal P}_{m,r}(|m|;\vec m;\vec c)$
($r\in\ZZ$)
with $-$ walls (resp. $+$ walls) and weight
$wt(b_0)=-m(\Lm_0-\Lm_1)-l\del$.
Then $b_0^*$ is an element of
${\cal P}_{m,l}$.
%%% because %%%%
%%% if $m>0$ and $-$ walls, $wt(b_0)=|m|(\Lm_1-\Lm_0)=-m(\Lm_0-\Lm_1)$
%%% if $m<0$ and $+$ walls, $wt(b_0)=|m|(\Lm_0-\Lm_1)=-m(\Lm_0-\Lm_1)$

This vector $b_0$ is in the following form;
$$
b_0=(\cd,0,0,d_1,d_2,\cd,d_{|m|-1},(-1)^{t+1}m,(-1)^{t+2}m,\cd),
\q
d_j=\underbrace{-j,j,-j,\cd}_{2c_j},
$$
where $t$ is the same one as in Theorem \ref{star-on-ext}.
Because of the form of $d_j$ and the fact $t(d_j)=j$,
the position of the left-most entry in any domain
is odd.
Then by Theorem \ref{star-on-ext},
we can write
$$
b_0^*=(\cd,0,0,d'_1,d'_2,\cd,d'_{|m|-1},(-1)^{t+1}m,(-1)^{t+2}m,\cd),
\q
d'_j=\underbrace{-j,j,-j,\cd}_{2c_j}.
$$
Thus, we have $d_j=d'_j$ and then
$b_0=b_0^*$ as an element of ${\cal P}_m(|m|;\vec m;\vec c)$.
(Note that in general, $wt(b_0)-wt(b_0^*)\in \ZZ\del$.).

\nd
Since ${\cal P}_{m,l}(|m|;\vec m;\vec c)$ is generated by
the extremal vector $b_0^*$ and $b_0^*\in \bmax(\lm)$, we get
$$
{\cal P}_{m,l}(|m|;\vec m;\vec c)\subset \bmax(\lm)
\q{\rm for\,\, any}\q
\vec c.
$$
Now, we get the desired result.\qed

\subsection{Description of $B(-\lm)^*$}

For $\lm=m(\Lm_0-\Lm_1)+l\del$ ($m,l\in\ZZ$),
the generator $u_{-\lm}=u_{\ify}\ot t_{-\lm}\ot u_{-\ify}$
corresponds to the path
$$
b_0=(\cd,\stackrel{-2}{0},\stackrel{-1}{0},\stackrel{0}{-m},
\stackrel{1}{m},\stackrel{2}{-m},\cd).
$$
We set
$$
 \vec 0:=\underbrace{(0,0,\cd,0)}_{|m|-1}.
$$
Then
$b_0$ is an element of
${\cal P}_{-m,-l}(|m|;-\vec m;\vec 0)$.
This implies
$$
B(-\lm)\cong {\cal P}_{-m,-l}(|m|;-\vec m;\vec 0).
$$
As we know by the results in the previous subsection,
$B(-\lm)\subset \bmax(-\lm)$. Thus, every element in $B(-\lm)^*$
is an extremal vector with weight $\lm$.
For an element $b\in {\cal P}_{-m,-l}(|m|;-\vec m;\vec 0)$,
let $d_1(b),d_2(b),\cd,d_{|m|-1}(b)$ be its finite domains and
$l_1(b),l_2(b),\cd,l_{|m|-1}(b)$ be their lengths.
Since any domain in ${\cal P}_{-m,-l}(|m|;-\vec m;\vec 0)$ is
a regular domain (see \cite[Sec 7.1]{N}),
recalling the definition of domain parameter (\cite[Sec 7.1]{N}), we have
$l(d_j)=0$ or 1 for any $j$, namely,
\beq
{\cal P}_{-m,-l}(|m|;-\vec m;\vec 0)
=\{b\in{\cal P}_{m,-l}(|m|)\,|\,
l(d_j)=0{\hbox{ or 1 for $j=1,2,\cd,|m|-1$}}\}.
\eeq
We can describe the action of $*$
on an element in $B(\lm)\subset \bmax(\lm)$ by Theorem \ref{star-on-ext}
because $*^{-1}=*$. Then we obtain the following;

\begin{pr5}
\label{b-lm-star}
For $\lm=m(\Lm_0-\Lm_1)+l\del$,
\beq
B(-\lm)^*
=\left\{b\in B(\til U_q(\slh)_0)|
\begin{array}{l}
{\hbox{ $b$ is an extremal vector with $|m|$ walls}},\,\,
wt(b)=\lm,\\
\q l_j(b)=0{\hbox{ or 1 for any $j=1,\cd,|m|-1$}}
\end{array}
\right\}.
\label{b-lm-star-formula}
\eeq
\end{pr5}

\subsection{Explicit form of Peter-Weyl type decomposition
of $B(\til U_q(\slh)_0)$}

The crystal of modified quantum algebra $B(\til U_q(\slh))$
has a decomposition of bi-crystal;
$$
B(\til U_q(\slh))=B(\til U_q(\slh)_+)\oplus
B(\til U_q(\slh)_0)\oplus B(\til U_q(\slh)_-),
$$
where $\til U_q(\slh)_{\pm}
:=\sum_{\pm\lan c,\lm\ran>0} U_q(\slh)a_{\lm}$.
The crystallized
Peter-Weyl type decomositios of $B(\til U_q(\slh)_{\pm})$
have been given in \cite{K4}.
By applying Proposition \ref{Peter-Weyl}, we can describe
the crystallized Peter-Weyl type decomposition for
$B(\til U_q(\slh)_0)$.

\begin{thm5}
\label{P-W-sl2}
There exists the following isomorphism of bi-crystal;
\beq
B(\til U_q(\slh)_0)\cong
\bigoplus_{\lm\in P_0/W}\bmax(\lm)\ot B(-\lm)^*,
\eeq
where $P_0=\{\lm\in P|\lan c,\lm\ran=0\}$ and $W$ is the Weyl group
associated with $\slh$.
\end{thm5}

\vskip6pt
\nd
{\sl Proof.\,}
In the course of the proof of Proposition \ref{Peter-Weyl},
it is enough to show that
the following conditions hold:
\begin{description}
\item[(C1')] For any extremal vector $b\in B(\til U_q(\slh)_0)$, there
exists an embeding of crystal
$$
B(wt(b))\hookrightarrow B(\til U_q(\slh)_0),
$$
given by $u_{wt(b)}\mapsto b$.
\item[(C2')] For any $\lm\in P_0$, $B(\lm)_{\lm}=\{u_{\lm}\}$.
\item[(C3')] For any extremal vectors $b_1,b_2\in B(\lm)$
($\lm\in P_0$) there exist $i_1,i_2,\cd,i_k$ such that
$$
b_2=S_{i_1}S_{i_2}\cd S_{i_k}b_1.
$$
\end{description}

\nd
(C1') \,\, For an extremal vector $b\in B(\til U_q(\slh)_0)$
let $B'$ be the connected component including $b$
and set $\lm=m(\Lm_0-\Lm_1)+l\del:=wt(b)$.
By Corollary 7.28 in \cite{N}, we have
\beq
\begin{array}{lllll}
B'&\cong& \aff(B^{\ot |m|})_{\bar l}&\cong& B(\lm)
\label{congcong}\\
b&\leftrightarrow& z^l\ot (\ep)^{\ot |m|}&\leftrightarrow& u_{\lm},
\end{array}
\eeq
where $\ep\in\{\pm\}$ and
$\bar l\in\{0,1,\cd,|m|-1\}$ and $\bar l\equiv l$ mod $|m|$.
This implies that (C1') holds.

\vskip5pt
\nd
(C2')
\,\,
The vector $u_{\lm}$ corresponds to the path
$$
b=(\cd,\stackrel{-2}{0},\stackrel{-1}{0},\stackrel{0}{m},
\stackrel{1}{-m},\stackrel{2}{m},\cd).
$$
This $b$ is an element of
${\cal P}_{m,l}(|m|;\vec m;\vec 0)$.
Thus, we have
$$
B(\lm)\cong
{\cal P}_{m,l}(|m|;\vec m;\vec 0).
$$
By Lemma 7.27 in \cite{N} and the comments below that lemma,
we get (C2').

\vskip5pt
\nd
(C3')\,\,
By the formula (7.47) in \cite{N}, we obtain the transitivity
of extremal vectors in $B(\lm)$ for $\lm \in P_0$. \qed

\begin{cor5}
Let  us denote $P_1(m,l;\vec c)$ for
${\cal P}_{m,l}(|m|;\vec m;\vec c)$ and
$P_2(m,l)$ for the R.H.S of (\ref{b-lm-star-formula}).
We have
$$
B(\wtil U_q(\slh)_0)
\cong
\bigoplus_{
m(\Lm_0-\Lm_1)+l\del\in P_0/W\atop
{\hbox{$\vec c$}}\in\ZZ_{\geq 0}^{|m|-1}}
P_1(m,l;\vec c)\ot P_2(m,l).
$$
\end{cor5}


\begin{thebibliography}{99}


\def\CMP{\sl Commum.Math.Phys.}
\def\IJMP{\sl Int.J.Mod.Phys.}
\def\Duke{\sl Duke Math.J.}

%\bibitem{BLM} Beilinson A A, Lusztig G. and MacPherson R,
%   A geometric setting for the quantum deformation of $GL_n$,
%             {\it Duke Math. J.}, {\bf 61}, (1990), 655--677.
%\bibitem{DFJMN} Davies B, Foda O, Jimbo M, Miwa T and Nakayashiki A,
%      Diagonalization of the XXZ Hamiltonian by vertex operators,
%              {\CMP}, {\bf 151}, (1993) 89--153.
%
%\bibitem{IIJMNT} Idzumi M, Iohara K, Jimbo M, Miwa T, Nakashima T and
%    Tokihiro T, Quantum affine symmetry in vertex models,
%    {\IJMP} A Vol 8,
%     No.8, (1993) 1479--1511.

%\bibitem{KMN}
%Kang S-J, Kashiwara M, Misra K, Miwa T, Nakashima T and Nakayashiki A,
%  Affine crystals and vertex models,
%  {\IJMP},{A7 }Suppl. 1A (1992) 449--484.
%
%\bibitem{Kac} Kac V G, Infinite dimensional Lie algebras,
%        3rd edition, Cambridge Univ. Press (1990).

\bibitem{K1} Kashiwara M,
 On crystal bases of the $q$-analogue of universal enveloping algebras,
        {\it Duke Math. J.},{\bf 63} (1991) 465--516.

\bibitem{K2} Kashiwara M,
             Global crystal bases of quantum groups,
             {\sl Duke Math. J.}, {\bf 69} (1993) 455--485.

\bibitem{K3} Kashiwara M,
      Crystal base and Littelmann's refined Demazure character formula,
             {\sl Duke Math. J.}, {\bf 71}(1993), 839--858.

\bibitem{K4} Kashiwara M,
      Crystal base of modified quantized enveloping algebra,
             {\sl Duke Math. J.}, {\bf 73} (1994), 383--413  .

\bibitem{K5} Kashiwara M,
        On Crystal Bases,
        Representations of Groups,
        Proceedings of a Summer Seminar held at Banff, Alberta,
        June 15 to 24, 1995
        (B.N. Allison and G.H. Cliff, eds), CMSAMS, Amer. Math. Soc.,
        Providence, RI (to appear).

%\bibitem{K6} Kashiwara M, private communication.

\bibitem{KKM} Kang S-J, Kashiwara M and Misra K.C,
        Crystal base of Verma modules for quantum Affine Lie algebras,
              {\sl Compositio Mathematica}, {\bf 92} (1994), 299--325.

\bibitem{KN} Kashiwara M and Nakashima T,
             Crystal graph for representations
             of the $q$-analogue of classical Lie algebras,
            {\sl J. Algebra}, {\bf 165}, (1994), 295--345.
\bibitem{L1} Lusztig G,
             Introduction to Quantum Groups,
             Birkh$\ddot{\hbox{a}}$user Boston, (1993).
\bibitem{L2} Lusztig G,
             Canonical bases in tensor product,
             {\sl Proc. Nat. Acad. Sci. USA}, {\bf 89}
             (1992), 8177--8179.
%\bibitem{T}  Nakashima T, Quantum R-matrix and Intertwiners for the
%      Kashiwara algebra,
%      {\sl Commun. Math. Phys.}, {\bf 164}, Number 2, (1994), 239-258.
\bibitem{N}  Nakashima T, Crystallized structure for level 0
      part of modified quantum affine algebra, (1995), q-alg 9506023,
      Osaka Univ.Math.Sci.preprint 4.

\end{thebibliography}
\end{document}